\documentclass[hidelinks,onefignum,onetabnum]{siamart251216}

\usepackage{lipsum}
\usepackage{amsfonts}
\usepackage{graphicx}
\usepackage{epstopdf}
\usepackage{algorithmic}
\usepackage{bm}
\usepackage{adjustbox}

\ifpdf
  \DeclareGraphicsExtensions{.eps,.pdf,.png,.jpg}
\else
  \DeclareGraphicsExtensions{.eps}
\fi

\newsiamremark{remark}{Remark}
\newsiamremark{hypothesis}{Hypothesis}
\crefname{hypothesis}{Hypothesis}{Hypotheses}
\newsiamthm{claim}{Claim}
\newsiamremark{fact}{Fact}
\crefname{fact}{Fact}{Facts}
\newsiamthm{problem}{Problem}

\headers{Discretization of the Surface Stokes-Helfrich model}{E. Igel, M. Porrmann, A. Voigt}

\title{A stream-function formulation for the Surface Stokes-Helfrich model - Surface Finite Element discretization and validation  \thanks{\funding{This work was funded by the German Research Foundation (DFG) within the research unit "Vector- and tensor-valued surface PDEs" - project numbers 417223351 (FOR3013).}}}

\author{Enno Igel\thanks{Institute of Scientific Computing, TUD Dresden University of Technology, Dresden, Germany
  (\email{enno.igel@tu-dresden.de}).}
\and Maik Porrmann\thanks{Institute of Scientific Computing, TUD Dresden University of Technology, Dresden, Germany
  (\email{maik.porrmann@tu-dresden.de}).}
\and Axel Voigt\thanks{Institute of Scientific Computing, TUD Dresden University of Technology, Dresden, Germany; Center for Systems Biology Dresden (CSBD), Dresden, Germany; Cluster of Excellence Physics of Life (PoL), TUD Dresden University of Technology, Dresden, Germany (\email{axel.voigt@tu-dresden.de}).}
}

\usepackage{amsopn}

\ifpdf
\hypersetup{
  pdftitle={A stream-function formulation for the Surface Stokes-Helfrich model - Surface Finite Element discretization and validation},
  pdfauthor={Enno Igel, Maik Porrmann, and Axel Voigt}
}
\fi
\newcommand{\R}{\mathbb{R}}
\renewcommand{\vec}[1]{\boldsymbol{#1}}
\newcommand{\normal}{\vec{\nu}}
\renewcommand{\S}{\mathcal{S}}
\renewcommand{\H}{\mathcal{H}}
\newcommand{\K}{\mathcal{K}}

\newcommand{\param}{\vec{X}}
\newcommand{\B}{\boldsymbol{B}}
\renewcommand{\P}{\boldsymbol{P}}
\newcommand{\Id}{\mathbf{Id}}
\newcommand{\tangent}{\mathrm{T}}

\newcommand{\flow}{\vec{u}}
\newcommand{\flowTest}{\vec{v}}

\newcommand{\pressure}{p}

\newcommand{\stress}{\vec{\sigma}}
\newcommand{\pos}{\vec{x}}

\newcommand{\Volume}{\Omega}
\newcommand{\tr}{\operatorname{tr}}

\newcommand{\gradS}{\nabla_{\!\!\S}}

\newcommand{\gradSVar}[1]{\nabla_{\!#1}}

\newcommand{\divS}{\operatorname{div_{\!\S}}}
\newcommand{\vecdivS}{\operatorname{\mathbf{div}_{\!\S}}}
\newcommand{\gradC}{\operatorname{\nabla_{\!C}}}
\newcommand{\vecgradC}{\operatorname{\bm\nabla_{\!C}}}

\newcommand{\divC}{\operatorname{div_{\!C}}}
\newcommand{\vecdivC}{\operatorname{\mathbf{div}_{\!C}}}

\newcommand{\DivC}{\operatorname{Div_{\!C}}}
\newcommand{\vecDivC}{\operatorname{\mathbf{Div}_{\!C}}}

\newcommand{\lapS}{\operatorname{\Delta_{\S}}}

\newcommand{\lapC}{\operatorname{\Delta_{C}}}
\newcommand{\veclapC}{{\bm\Delta_{C}}}

\newcommand{\curlS}{\operatorname{curl_{\mkern-1mu\S}}}
\newcommand{\veccurlS}{\mathbf{curl}_{\mkern -1mu \S}}
\newcommand{\veccurlSVar}[1]{\mathbf{\operatorname{curl_{\mkern-1mu#1}}}}

\newcommand{\rateOfDeformation}{\bm \sigma}

\newcommand{\areaelement}[1]{d\S \ifthenelse{\equal{#1}{}}{}{_{#1}}}

\newcommand{\Ltwo}{\mathrm{L}^2}

\newcommand{\TangVectors}{\tangent\S}
\newcommand{\Vectors}{\restrict{\tangent\R^3}{\S}}

\newcommand{\nTangTensors}[1]{\tangent^{#1}\S}
\newcommand{\Scalars}{\nTangTensors{0}}

\newcommand{\nTensors}[1]{\restrict{\tangent^{#1}\R^3}{\S}}

\newcommand{\Inner}[3][]{\left({#2}\,,\,{#3}\right)\ifthenelse{\equal{#1}{}}{}{_{#1}}}

\newcommand{\restrict}[2]{\left.#1\right\vert_{#2}}

\newcommand{\reducedVolume}{V_r}

\usepackage{mathtools}
\newcommand{\definedAs}{\coloneqq}

\newcommand{\amdis}{\texttt{AMDiS}}
\newcommand{\dune}{\textsc{Dune}}
\newcommand{\alugrid}{\textsc{Dune-alugrid}}
\newcommand{\curvedgrid}{\textsc{Dune-curvedgrid}}
\newcommand{\functions}{\textsc{Dune-functions}}
\newcommand{\eigen}{Eigen}

\begin{document}

\maketitle

\begin{abstract}
We introduce a stream-function formulation for the Surface Stokes-Helfrich model and study an Arbitrary Lagrangian-Eulerian (ALE) Surface Finite Element (SFEM) discretization in space and a semi-implicit discretization in time of these equations. In contrast with special cases of a surface Stokes model on stationary or prescribed evolving surfaces the surface pressure remains as an unknown and needs to be reconstructed in each time step. As for the surface Stokes model on a prescribed evolving surface the formulation also requires an additional gradient potential. A detailed computational comparison with an ALE-SFEM discretization in space and a similar semi-implicit discretization in time of the corresponding velocity-pressure formulation is explored on a representative problem formulation for simply-connected surfaces, demonstrating the validity of the approach.
\end{abstract}



\section{Introduction}

We address fluid deformable surfaces, which are soft materials exhibiting a solid–fluid duality. They store elastic energy when stretched or bent, like solid shells, but cannot do so under in-plane shear, a situation under which they flow as two-dimensional, viscous fluids. This duality is responsible for a tight coupling between tangent flows and shape changes in the presence of curvature and thus makes curvature a natural element of such materials, which can be found in living objects at length scales ranging from cell membranes to single cells to tissues and organs \cite{schamberger2023curvature}. The basic model describing this behavior is a Surface (Navier-)Stokes-Helfrich model \cite{arroyo2009relaxation,torres2019modelling,reuther2020numerical,nitschke2026hydrodynamic}. It combines bending and surface hydrodynamics and assumes the surface fluid to be inextensible. The unknowns are the surface velocity (in tangential and normal direction) and the surface pressure. The model has been mainly applied to simulate the shape evolution of lipid membranes \cite{torres2019modelling,reuther2020numerical,bachini2023derivation,sischka2025two,sahu2025arbitrary,nitschke2026hydrodynamic}. As a result of the small length scale and the relatively large viscosity of lipid membranes a Stokes approximation is in most cases sufficient for this application.

As in previous studies we investigate the most common situation of closed surfaces and as in \cite{krause2023numerical,nestler2023stability,KV24,PV24,sischka2025two,garcke2025parametric,PBV25} we also consider a constraint on the enclosed volume, which makes the problem nonlocal. This model, in its full generality or under various symmetry assumptions, has been considered numerically in~\cite{torres2019modelling,reuther2020numerical,al2020shear,krause2023numerical,olshanskii2023equilibrium,sauer2025curvilinear,garcke2025parametric}. For the general setting the computational approaches are primarily based on an Arbitrary Lagrangian-Eulerian (ALE) approach \cite{reuther2020numerical,sauer2025curvilinear,sahu2025arbitrary}, an isoparametric Surface Finite Element Method (SFEM) \cite{DE13,nestler2019finite}, a method to move grid points in tangential direction to maintain mesh regularity \cite{Barrett_SIAMJSC_2008,sauer2025curvilinear,PBV25}, a stable element pair for the surface velocity and surface pressure, e.g. a Taylor-Hood element, and a semi-implicit time discretization. Stability estimates could be shown for such schemes in \cite{garcke2025parametric}. However, detailed analysis of the schemes is not yet available. This drastically changes if the surface is not evolving or if its evolution is prescribed. In these situations the problem changes to a Surface (Navier-)Stokes equation modeling a two-dimensional fluid with the only unknowns being the tangential velocity and the surface pressure. For related numerical investigations we refer to \cite{reuther2015interplay,reuther2018erratum,Fries,lederer2020divergence,bonito2020divergence,brandner2022finite,neilan2025ac,elliott2025sfem}
for stationary surfaces and \cite{olshanskii2022tangential,nestler2023stability,olshanskii2023eulerian,elliott2025fully} for prescribed evolving surfaces. For the former case the formulation is either based on the classical velocity-pressure formulation or on a reformulation using vorticity and stream-function as unknowns. Similar to the classical two-dimensional incompressible Navier-Stokes equations the stream-function formulation circumvents the saddle-point structure and allows to solve the problem by solving a higher order scalar-valued partial differential equation. Due to an increased complexity for vector-valued partial differential equations on surfaces \cite{nestler2019finite,hardering2023tangential,bachini2024diffusion}, the latter point is even more relevant on surfaces. Various attempts are based on such a reformulation \cite{10.1145/1189762.1189766,Nitschke_Voigt_Wensch_2012,reuther2015interplay,reuther2018erratum,10.1093/imanum/dry062,GROSS2018663,torres2019modelling,pearce2019geometrical}.
However, they should be considered with care. The approach is based on a decomposition of the velocity field in divergence- and
curl-free parts. While appropriate on simply connected surfaces, this is not sufficient on non-simply connected surfaces,
e.g. a torus. As a consequence of the topology also non-trivial harmonic parts, velocity fields which are divergence- and curl-free, exist, see \cite{nitschke2021vorticity} for a detailed discussion. Only recently the approach has been extended to general surfaces by solving also for the harmonic vector field \cite{FluidCohomology,ZhuYinChern2025,bruers2025streamfunction}. While computationally more demanding, the formulation remains algorithmically attractive, as it still avoids a saddle-point structure and only involves primarily scalar unknowns. Extensions of the stream-function formulation towards evolving surfaces, which go beyond prescribed surface evolutions, as considered in  \cite{reuther2015interplay,reuther2018erratum}, do not exist, not even for simply connected surfaces. We here explore this situation and compare the results with existing velocity-pressure formulations for the Surface Stokes-Helfrich model.

The paper is organized as follows. We introduce the notation and provide the various formulations for the Surface Stokes-Helfrich model in
\cref{sec:model}, formulate the numerical approach in \cref{sec:numerics}, compare the results with a velocity-pressure formulation on a benchmark problem in \cref{sec:results}, and draw conclusions in \cref{sec:conclusions}. A detailed derivation of the stream-function formulation is provided in the appendix \ref{app:derivation}.

\section{Modeling}
\label{sec:model}

We use a similar notation as in \cite{bachini2023derivation,PBV25}, which is here repeated for convenience. We consider a time dependent smooth and oriented surface $\S = \S(t)$ without boundary, embedded in $\R^3$ and given via a parametrization $\param = \param(t)$. The enclosed volume is denoted by $\Volume = \Volume(t)$. We denote by $\normal$ the outward pointing surface normal, the surface projection is $\P=\Id -\normal \otimes \normal$, with $\Id$ the identity matrix, the shape operator is $\B= -\vecgradC \normal$, the mean curvature $\H = \tr \B$, and the Gaussian curvature $\K = \frac{1}{2}\left(\H^2-\|\B\|^2\right)$. We consider time-dependent Euclidean-based $ n $-tensor fields on $\S$ denoted by $\nTensors{n}$. We call $\nTensors{0} = \Scalars$ the space of scalar fields, $\nTensors{1}=\Vectors$ the space of vector fields, and $\restrict{\tangent^2\R^3}{\S} $ the space of 2-tensor fields on $\S$. Important subtensor fields are tangential $n$-tensor fields in $\nTangTensors{n} \subset \nTensors{n}$.\\
Let $p= p(\pos,t) \in \tangent\S(t)$ be a continuously differentiable scalar field, $\flow = \flow(\pos, t) \in \restrict{\tangent\R^3}{\S(t)}$ a continuously differentiable $\R^3$-vector field, and $\stress = \stress(\pos, t) \in \restrict{\tangent^{2}\R^3}{\S(t)}$ a continuously differentiable $\R^{3\times3}$-tensor field, each defined on $\S(t)$. We define the (componentwise) surface gradient by $\gradC p = \nabla p^e\P$, $\vecgradC\flow = \bm\nabla\flow^e\P$ and $\vecgradC\stress = \bm\nabla\stress^e \P$,
where $p^e$, $\flow^e$ and $\stress^e$ are arbitrary smooth extensions of $p$, $\flow$ and $\stress$ in the normal direction and $\nabla$ and $\bm\nabla$ are the gradients of the embedding space $\R^3$. The traces of those operators give rise to one notion of divergence operators, which for a vector field $\flow$ and a $2$-tensor field $\stress$ are $\DivC\flow = \tr(\vecgradC\flow)$ and $\vecDivC \stress = \tr_{(2,3)}(\vecgradC\stress)$, where $\tr$ is the trace operator.
However, those divergences are not the adjoints of the gradients, which we will denote by $\divC = -\gradC^*$ if acting on vectors and $\vecdivC = -\vecgradC^*$ if acting on $2$-tensors.
The relations to the covariant derivative $\gradS$ and the covariant divergence $\divS$ on $\S$, read $\gradC p=\gradS p$ and ${\DivC\flow = \divS(\flow_T)-u_N\H}$, in contrast to $\divC \flow = \divS( \flow_T)$, with $\flow_T = \P\flow$ and $u_N = \flow \cdot \normal$. We only use covariant curl operators, namely $\veccurlS : \Scalars \to \TangVectors, \veccurlS p \definedAs \normal \times \gradS p$ and its negative adjoint $\curlS : \TangVectors \to \Scalars, \curlS (\flow_T) \definedAs \divS(\flow_T \times \normal)$.  Lastly we can define the componentwise Laplace operators $\lapC = \divC\circ\gradC$ and $\veclapC = \vecdivC\circ\vecgradC$. Note that only $\lapC$ agrees with the Laplace-Beltrami operator $\lapS = \divS \circ \gradS$, while $\veclapC$ is used to write the $(2,3)$-contracted Gauss-Weingarten equations as $\H\normal=\veclapC \param$ with a common abuse of notation. With these notations we are now ready to formulate the models of interest.

The Surface Stokes-Helfrich model in velocity-pressure formulation reads \cite{torres2019modelling,PBV25}:
\begin{problem}
\label{prob:u-p}
Given an initial (closed, simply connected) surface $\S(0)$ via $\param(0)$, determine the parametrization $\param$, the surface velocity $\flow$, the surface pressure $p$ and the Lagrange multiplier $\lambda$, such that
\begin{subequations}
\begin{align}
    \gradS p + p \H \normal -  \mu\vecdivC \rateOfDeformation(\flow) + \lambda \normal + \gamma \flow&= f_B \normal  \label{eq:u-p-mb} \\
    \DivC \flow &= 0 \label{eq:u-p-inext}\\
    \int_\S \flow\cdot \normal \ \areaelement{} &= 0 \label{eq:volume_conservation}
\end{align}
on $\S(t)$ and $\partial_t \param = (\flow \cdot \normal) \normal$, with bending force
\begin{align} \label{eq:bendingforce}
    f_B \normal &= - \kappa (\lapS \H + \frac{1}{2}\H^3 - 2 \H \K) \normal \,.
\end{align}
\end{subequations}
Here, $\rateOfDeformation(\flow) \definedAs (\P\vecgradC\flow + \left(\P\vecgradC\flow\right)^T )$ is the surface rate of deformation tensor, $\mu > 0$ denotes the surface viscosity, $\gamma > 0$ is the friction coefficient, and $\kappa > 0$ the bending rigidity.
\end{problem}

In Problem~\ref{prob:u-p}, \eqref{eq:u-p-mb} is the force balance, while \eqref{eq:u-p-inext}, together with the Lagrange multiplier $\pressure$, enforces the inextensibility constraint. Lastly, \eqref{eq:volume_conservation} enforces the enclosed volume to remain constant and $\lambda$ is the corresponding Lagrange multiplier. Due to the reparametrization invariance of the surface, we only use the normal part $(\flow \cdot \normal)\normal$ of the material velocity to update the parametrization $\param$.

Considering a tangential-normal splitting of the surface velocity $\flow = \flow_T + u_N \normal$, where $\flow_T = \P \flow$ is the tangetial part and $u_N = \flow \cdot \normal$ the normal component of $\flow$, allows one to rewrite the equations. The new formulation reads \cite{reuther2015interplay,reuther2018erratum,jankuhn2018incompressible,reuther2020numerical}:
\begin{problem}
    \label{prob:tang-nor-split}
Given an initial (closed, simply connected) surface $\S(0)$ via $\param(0)$, determine the parametrization $\param$, the surface tangential velocity $\flow_T$, the surface normal velocity $u_N$, the surface pressure $p$ and the Lagrange multiplier $\lambda$, such that
\begin{subequations}
\begin{align}
    \gradS p - \mu \left[ \vecdivS \rateOfDeformation(\flow_T) - 2  \vecdivS(u_N \B)\right] + \gamma \flow_T&= \vec{0}  \label{eq:tmb}\\
     p\H - \mu \left[\normal \cdot \vecdivC \rateOfDeformation(\flow_T) - 2 \normal \cdot \vecdivC(u_N \B)\right]  + \lambda + \gamma u_N &= f_B \label{eq:nmb}\\
    \divS \flow_T - u_N \H &= 0 \label{eq:incomp} \\
    \int_\S u_N \areaelement{} &= 0 \,,
\end{align}
\end{subequations}
on $\S(t)$ and $\partial_t \param = u_N \normal$, with $f_B$ defined as in \cref{eq:bendingforce}.
 \end{problem}

In Problem \ref{prob:tang-nor-split} the force balance has been split into \eqref{eq:tmb} and \eqref{eq:nmb}, while the inextensibility and volume constraints remain as before, only changing in notation. In this formulation the tight coupling between the surface tangential velocity $\flow_T$ and shape changes (surface normal velocity $u_N$) in the presence of curvature is revealed more clearly. This formulation is also the starting point to obtain the stream-function formulation.

To do so, we consider the Helmholtz decomposition for the tangential surface velocity $\flow_T \in \bm L_T^2(\S) = \{\bm f \in \bm L^2(\S), \bm f \cdot \normal = 0 \}$ with $\bm L^2(\S) = [L^2(\S)]^3$, which can be uniquely decomposed into
\begin{align}
    \flow_T = \veccurlS \phi + \gradS \psi + \bm \xi \label{eq:Helmholtz}
\end{align}
with $\phi \in H_*^1(\S)$ the stream-function, $\psi \in H_*^1(\S)$ the gradient potential, and $\bm \xi \in \bm L_{T,\text{har}}^2(\S)$ the harmonic part. Thereby $H_*^1 = \{ f \in H^1, \int_\S f \, \areaelement{} = 0\}$ and $\bm L_{T,\text{har}}^2(\S) = \{\bm f \in \bm L_T^2(\S), \divS \bm f = \curlS \bm f = 0\}$. In our case of a simply connected surface the harmonic part vanishes.

Following \cite{reuther2015interplay,reuther2018erratum}, we insert this decomposition into the mass conservation law \eqref{eq:incomp}, which gives $\lapS \psi  - u_N \H = 0$, an equation for $\psi$. The equations for the stream-function $\phi$ and the vorticity $\omega = \lapS \phi$ then essentially result by applying the $\curlS$ operator to \eqref{eq:tmb}. The equation for the surface normal velocity $u_N$ follows by inserting the Helmholtz decomposition in \eqref{eq:nmb}. The missing equation for the surface pressure $p$ essentially follows by applying the $\divS$ operator to  \eqref{eq:tmb}. The resulting equation can be seen as a pressure reconstruction, extending the approach in \cite{reuskenSfAnalysis,brandner2022finite}. All steps require intensive use of various identities. Details are provided in the appendix \ref{app:derivation}. The obtained formulation reads:

\begin{problem} \label{problem:sf}
Given an initial (closed, simply connected) surface $\S(0)$ via $\param(0)$, determine the parametrization $\param$, the surface vorticity $\omega$, the surface stream-function $\phi$, the surface gradient potential $\psi$, the surface normal velocity $u_N$, the surface pressure $p$, and the Lagrange multiplier $\lambda$, such that
\begin{subequations}
\allowdisplaybreaks
\begin{align}
    \mu\lapS \omega +2\mu\divS(\K \gradS \phi) + 2\mu \curlS(\K \gradS \psi)
    \nonumber \\
    - 2\mu\curlS(\B \gradS u_N + u_N \gradS \H) - \gamma \omega &= 0 \label{eq:sf_1}
    \\ 
    2 \mu\lapS(u_N \H) + 2\mu\divS(\K\gradS \psi) + 2\mu \divS(\K \veccurlS \phi)
    \nonumber \\
    - 2\mu \divS(\B \gradS u_N + u_N \gradS \H)
    - \lapS p - \gamma u_N \H &= 0 \label{eq:sf_2} \\ 
    -2 \mu \divS(\B (\gradS \psi + \veccurlS \phi)) + 2 \mu (\gradS \psi + \veccurlS \phi) \cdot \gradS \H \nonumber \\
    2 \mu u_N \B : \B + p\H + \gamma u_N +\lambda - f_B  &= 0 \label{eq:sf_3}
    \\ 
    \lapS \phi - \omega&= 0
    \label{eq:sf_4} \\ 
    \lapS \psi - u_N \H &= 0 \label{eq:sf_5} \\ 
    \int_\S u_N \ \areaelement{} &= 0 \label{eq:sf_6} \, , 
\end{align}
\end{subequations}
on $\S(t)$ and $\partial_t \param = u_N \normal$, where again $f_B$ is defined via \eqref{eq:bendingforce}.
\end{problem}

In Problem \ref{problem:sf} the force balance is expressed in \eqref{eq:sf_1} - \eqref{eq:sf_3}. \eqref{eq:sf_4} is an auxiliary equation, added in order to avoid higher order operators, the inextensibility becomes an equation for the gradient potential and the volume constraint remains as before. This formulation significantly differs from the known stream-function formulations where the surface is not evolving or where the evolution is prescribed. The main difference is the presence of the surface pressure $p$ in \eqref{eq:sf_3}. However, if we restrict the equations to these situations, the known formulations are recovered. For the case of a stationary surface ($u_N = 0$), \eqref{eq:sf_1} and \eqref{eq:sf_4} reduce to the system of equations for the surface vorticity $\omega$ and surface stream-function $\phi$ considered in \cite{reuskenSfAnalysis,brandner2022finite}. The gradient potential $\psi$ vanishes and \eqref{eq:sf_2} reduces to the pressure reconstruction considered in \cite{reuskenSfAnalysis,brandner2022finite}. For a prescribed evolution of the surface ($u_N$ given), \eqref{eq:sf_1}, \eqref{eq:sf_4} and \eqref{eq:sf_5} reduce to the system of equations for the surface vorticity $\omega$, surface stream-function $\phi$ and gradient potential $\psi$ considered in \cite{reuther2018erratum} if inertia terms are dropped in that formulation. Also in this case \eqref{eq:sf_2} can be seen as a pressure reconstruction. In both special cases the surface pressure $p$ is not required to solve the system. In contrast, the full problem requires $p$ to evolve the surface normal velocity $u_N$ in \eqref{eq:sf_3}.

We summarize that, even if the system of equations looks more complicated with stronger couplings with geometric quantities and requires additional unknowns, some advantages of the classical stream-function formulation remain. In particular, all unknowns are scalar-valued and the saddle-point structure for the pressure has been eliminated.

\section{Numerics}
\label{sec:numerics}
We exploit the reparametrization invariance of the surface numerically and consider an Arbitrary Lagrangian-Eulerian (ALE) Surface Finite Element Method (SFEM) \cite{DE13} to solve the highly nonlinear set of geometric and surface partial differential equations \eqref{eq:sf_1}-\eqref{eq:sf_6} and \eqref{eq:bendingforce} by combining the system with a mesh redistribution approach, see \cite{Barrett_SIAMJSC_2008}. Following the approaches in \cite{bachini2023derivation,krause2023numerical,PBV25}, the additional equations for the parametrization evolution are
\begin{subequations}
  \begin{align}
    \partial_t \param \cdot \normal &= u_N \label{eq:normalvel}\\
    \H \normal &= \veclapC \param \,. \label{eq:meandiff}
  \end{align}
\end{subequations}
  These equations not only address the movement in normal direction, but also generate a tangential mesh movement to maintain the shape regularity, and additionally provide an implicit representation of the mean curvature $\H$. We consider a discrete $k$-th order approximation $\S_h^k$ of $\S$, with $h$ denoting the discretization in space. We use the \curvedgrid{} library \cite{praetorius2020dunecurvedgrid} and consider each geometrical quantity like the normal vector $\normal_h$, the shape operator $\B_h$, the Gaussian curvature $\K_h$, and the $\Ltwo$-inner products $\langle\cdot , \cdot\rangle_h$ with respect to $\S_h^k$. In the following we will drop the index $k$ and write $\S_h$. We define the discrete function spaces for scalar functions by $V_{k}(\S_h)=\{ \psi \in C^0(\S_h) \vert\, \psi\vert_{T}\in\mathcal{P}_{k}(T)\}$ and for vector fields by $\boldsymbol{V}_{k}(\S_h)=[V_{k}(\S_h)]^3$. Within these definitions $T$ is the mesh element and $\mathcal{P}_{k}$ are the polynomials of order $k$. We consider $\param_h\in\boldsymbol{V}_3(\S_h)$ and  $\omega_h, \phi_h, \psi_h, u_{N,h}, p_h, \H_h \in V_3(\S_h)$. This leads to an isoparametric setting, suggested to be optimal for the surface Stokes equation on a stationary surface \cite{reuskenSfAnalysis,hardering2025parametric} and is consistent with the velocity-pressure formulation in \cite{PBV25}, which will be used for comparison.

  For the discretization in time we consider a first order implicit-explicit splitting resulting in a linear system that is solved at each time step. We treat linear terms implicitly, whereas non-linear terms are split such that geometric properties are evaluated at the previous time step, or in case of the mean curvature $\H_h$ are treated semi-implicitly. Given some discrete time steps $\{t_n\}_{n=0}^M,\,t_0=t_\text{start}, \, t_M=t_\text{end}$, we define $\bm Y_h^{n+1} \definedAs \param_h^{n+1} - \param_h^n$ and $\tau_{n+1} \definedAs t_{n+1} - t_{n}$.

  The resulting time-discrete linear system can be formulated for an initial surface $\param_h^0$ and corresponding discrete geometric properties $\normal_h^0, \B_h^0,\H_h^0, \K_h^0$ as well as $u_{N,h}^0$ and reads in each time step:
  \begin{problem}
  Given a discrete surface $\S_h^n$ via $\param_h^n$ with discrete geometric properties $\normal_h^n, \B_h^n,\H_h^n, \K_h^n$, determine $\omega_h^{n+1}, \phi_h^{n+1}, \psi_h^{n+1}, u_{N,h}^{n+1}, p_h^{n+1}, \H_h^{n+1} \in V_3(\S_h^n)$, $\lambda^{n+1}$ $\in \R$, and $\bm Y_h^{n+1} \in \bm V_3(\S_h^n)$, such that
\allowdisplaybreaks
\begin{align*}
     -\left< \mu \gradS \omega_h^{n+1},\, \gradS \widehat\phi \right>_h - \left< (2\mu \K_h^{n}-\gamma) \gradS \phi_h^{n+1},\, \gradS \widehat\phi \right>_h\\ - \left<2\mu \K_h^n \gradS \psi_h^{n+1} ,\, \veccurlS \widehat{\phi} \right>_h
    + \left< 2\mu (\B_h^n \gradS u_{N,h}^{n+1} + u_{N,h}^{n+1} \gradS \H_h^n ),\, \veccurlS \widehat{\phi} \right>_h &= 0 \\
    -\left< 2\mu \gradS (u_{N,h}^{n+1} \H_h^{n}) ,\, \gradS \widehat \psi\right>_h - \left< (2\mu \K_h^n - \gamma) \gradS \psi_h^{n+1},\, \gradS \widehat \psi \right>_h \\
    - \left<2\mu \K_h^{n} \veccurlS \phi_h^{n+1} ,\, \gradS \widehat{\psi} \right>_h
    + \left< 2\mu (\B_h^n \gradS u_{N,h}^{n+1} + u_{N,h}^{n+1} \gradS \H_h^n ),\, \gradS \widehat{\psi} \right>_h \\ + \left< \gradS p_h^{n+1} ,\, \gradS \widehat \psi\right>_h &= 0 \\
    \left<\gradS \psi_h^{n+1} + \veccurlS \phi_h^{n+1},\, 2\mu (\B_h^n \gradS \widehat{u_N} + \widehat{u_N} \gradS \H_h^n) \right>_h \\
    + \left< 2\mu u_{N,h}^{n+1} \B_h^{n},\, \widehat{u_N} \B_h^{n} \right>_h + \left< p_h^{n+1}\H_h^{n},\, \widehat{u_N}\right>_h + \left< \lambda^{n+1} ,\,\widehat{u_N} \right>_h \\
    - \kappa \left< \gradS \H_h^{n+1} ,\, \gradS \widehat{u_N} \right>_h + \kappa \left< \tfrac{1}{2}\H_h^{n+1}(\H_h^n)^2 - 2H_h^{n+1} \K_h^n ,\, \widehat{u_N} \right>_h&= 0 \\
    \left< \gradS\phi_h^{n+1},\, \gradS \widehat \omega\right>_h + \left< \omega_h^{n+1} ,\, \widehat\omega \right>_h &= 0 \\
    \left< \gradS \psi_h^{n+1},\, \gradS \widehat p \right>_h + \left< u_{N,h}^{n+1} \H_h^{n} ,\, \widehat p \right>_h + \left< u_{N,h}^{n} \H_h^{n+1} ,\, \widehat p \right>_h - \left< u_{N,h}^{n} \H_h^{n} ,\, \widehat p \right>_h&=0 \\
    \left< \bm Y_h^{n+1} \cdot \normal_h^n,\,\widehat{\H} \right>_h  - \tau_{n+1} \left< u_{N,h}^{n+1},\, \widehat{\H}\right>_h &= 0 \\
    \left< \H_h^{n+1} \normal_h^n,\, \widehat{\param} \right>_h + \left< \gradC \bm Y_h^{n+1},\, \gradC \widehat{\param} \right>_h + \left< \gradC \param^{n},\, \gradC \widehat{\param} \right>_h&= 0 \\
    \left< u_{N,h}^{n+1},\, \widehat \lambda\right>_h &= 0
\end{align*}
for all test functions $\widehat\omega, \widehat\phi, \widehat\psi, \widehat{u_{N}}, \widehat p, \widehat{\H} \in V_3(\S_h^n),$
$ \widehat \lambda \in \R, \widehat{\bm Y} \in \bm V_3(\S_h^n)$.
  \end{problem}

All differential operators are evaluated with respect to the discrete surface at time step $t_n$. For readability, we simply write $\gradS$ and  $\veccurlS$, which should be understood as $\gradSVar{\S_h^{n}}$ and $\veccurlSVar{\S_h^{n}}$, respectively. The surface is then updated by $\param^{n+1}_h = \param^{n}_h + \bm Y^{n+1}_h $ and all unknowns are lifted to the new surface. In order to recover the properties of $\phi$ and $\psi$ from the Helmholtz decomposition we project the solutions $\phi_h^{n+1}, \psi_h^{n+1}$ onto $V_{3,*}(\S_h^{n+1})$ by subtracting the mean. Afterwards, we reconstruct the tangential surface velocity for visualization and comparison:
\begin{problem}
Find $\flow_{T,h}^{n+1} \in \bm V_3(\S_h^{n+1})$, such that
\[
\langle \flow_{T,h}^{n+1}, \widehat{\flow_T} \rangle_h = \langle (\gradS \psi_h^{n+1} + \veccurlS \phi_h^{n+1}), \widehat{\flow_T} \rangle_h
\]
for all test functions $\widehat{\flow_T} \in \bm V_3(\S_h^{n+1})$.
\end{problem}

The velocity-pressure formulation is discretized as in \cite{PBV25}, using a $(\mathcal{P}_3,\mathcal{P}_2)$ Taylor-Hood element for the velocity-pressure pair and an isoparametric approach ($k = 3$) for the remaining unknowns, a similar approach to maintain mesh regularity and a similar implicit-explicit time discretization. For details and validation we refer to \cite{PBV25}.

We have realized the described discretizations within the finite element toolbox \amdis{}  \cite{VeyVoigt2007AMDiS,AmdisOld}, available at \cite{AMDiS:2.10}, which is a high level module in the \dune{} \cite{Dune2.10} ecosystem. For an overview of \dune{} we refer to \cite{duneBook}. Moreover, we use the modules \alugrid{} \cite{dune-alugrid} to manage a grid with sphere topology, \curvedgrid{} \cite{dune-curvedGrid} to map this grid onto the surface by the parametrization $\param_h$. This parametrization, as well as all other discrete functions and their bases are represented by means of \functions{} \cite{dune-functions}. We use \eigen{} \cite{eigenweb} as a linear algebra backend, and the direct solver PARDISO. Finally, all 3D visualizations have been created with ParaView \cite{Paraview}.

\section{Results}
\label{sec:results}

\begin{figure}[htpb]
\centering

\begin{adjustbox}{width=0.9\textwidth}
		\includegraphics{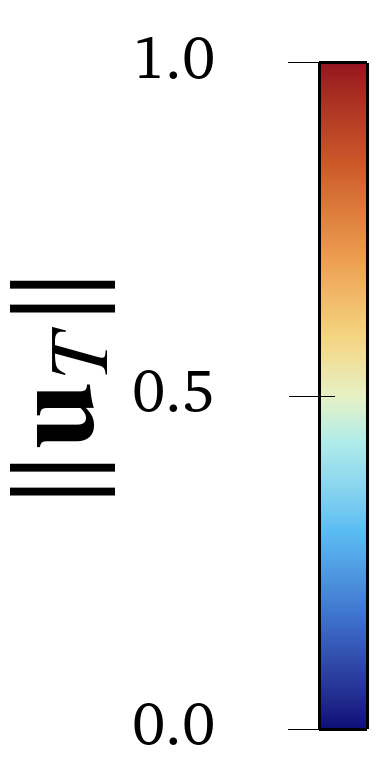}
		\includegraphics{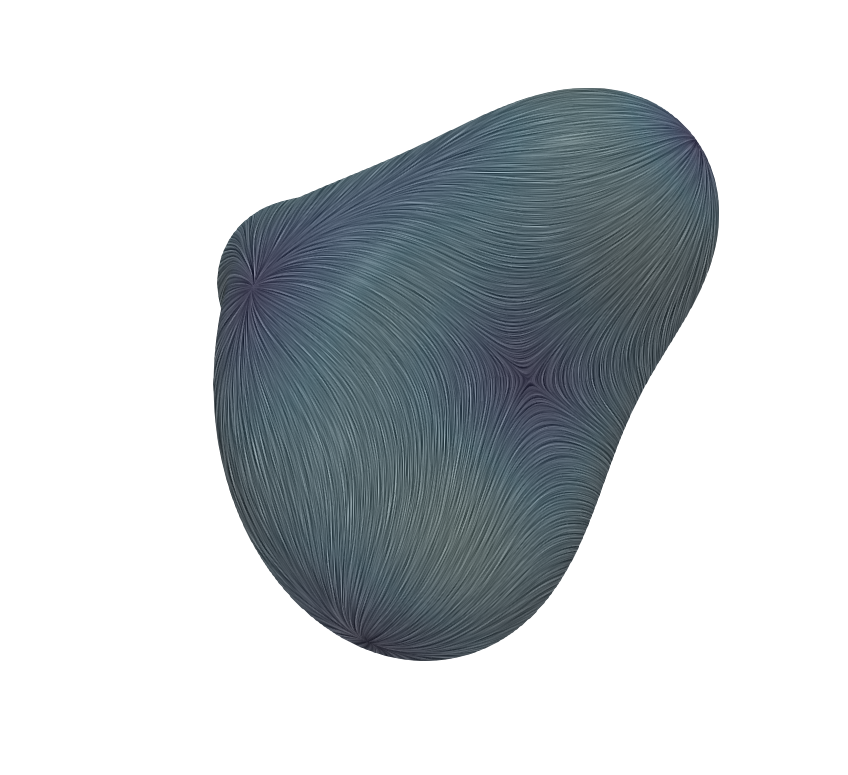}
		\includegraphics{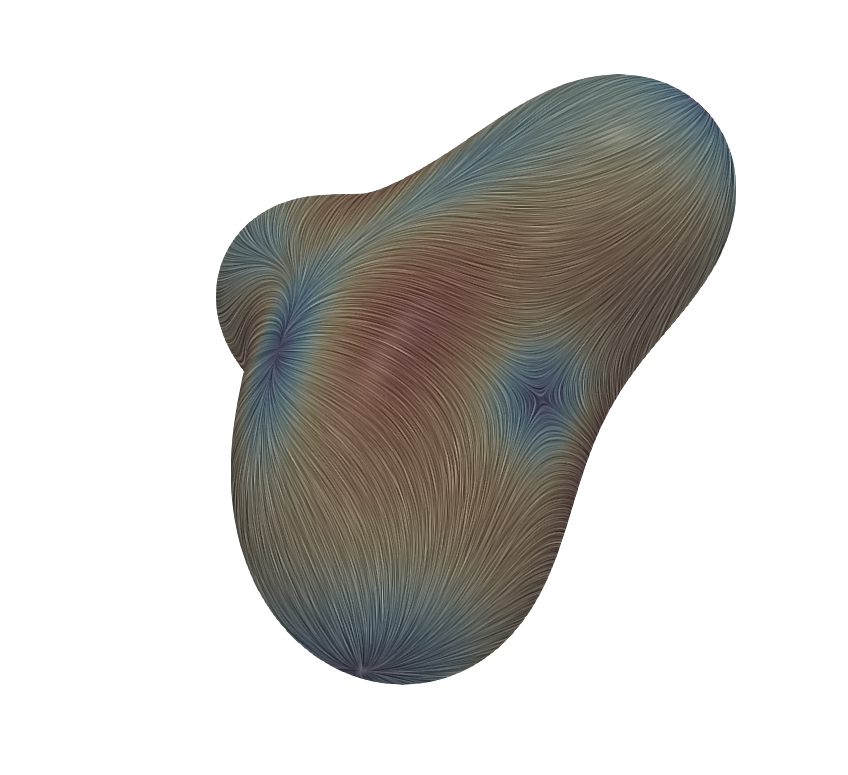}
		\includegraphics{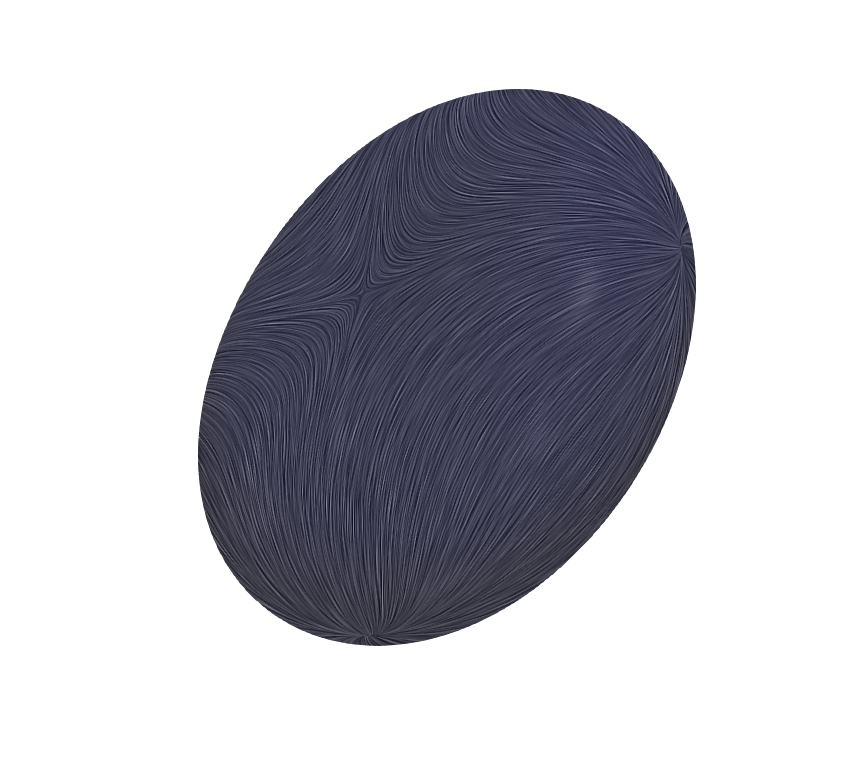}
		\includegraphics{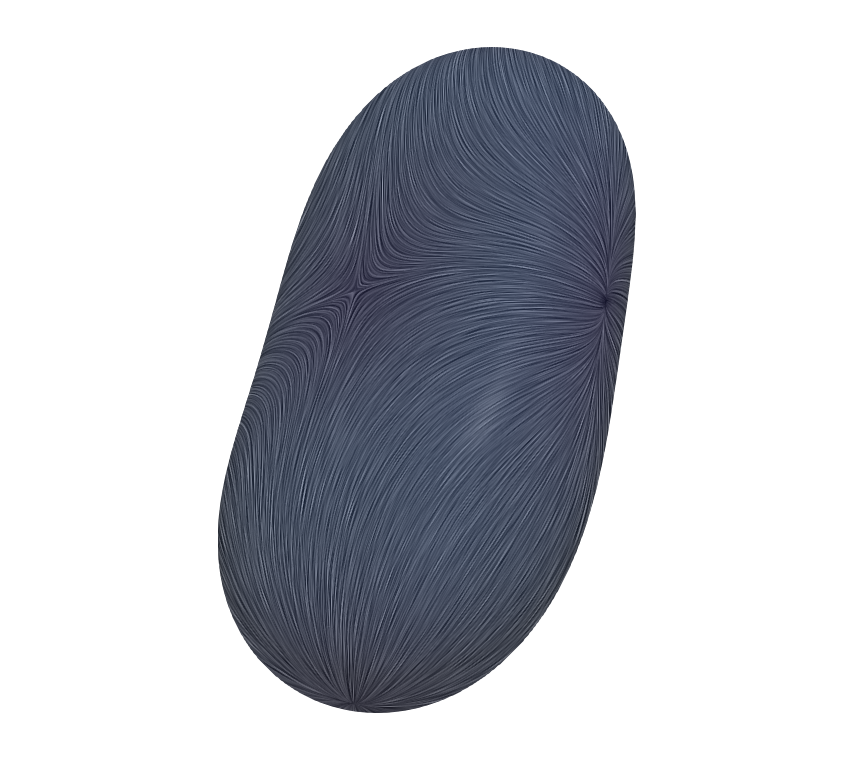}
		\includegraphics{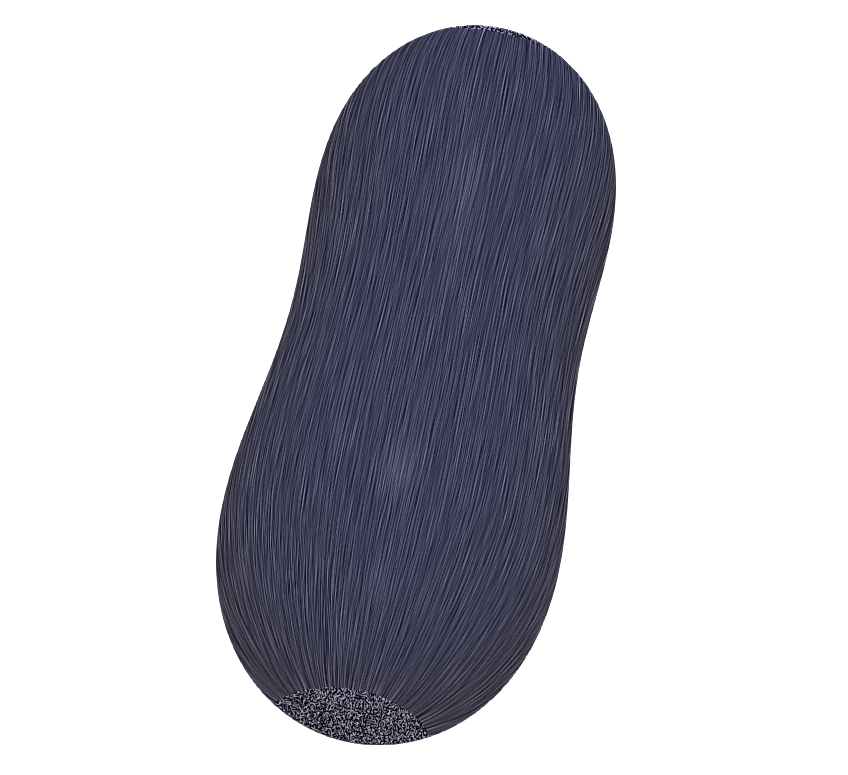}
	\end{adjustbox}

	\begin{adjustbox}{width=0.9\textwidth}
		\includegraphics{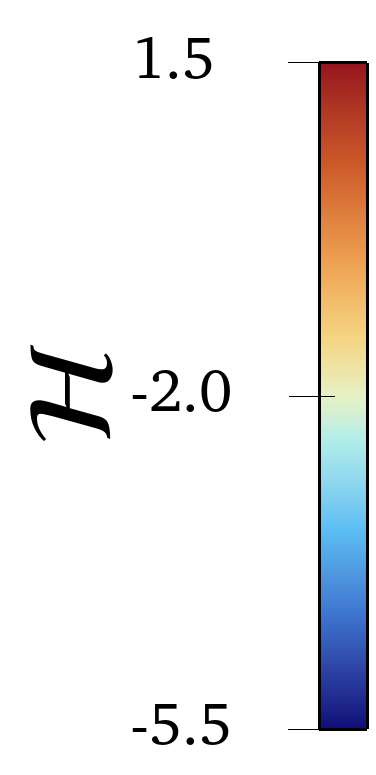}
		\includegraphics{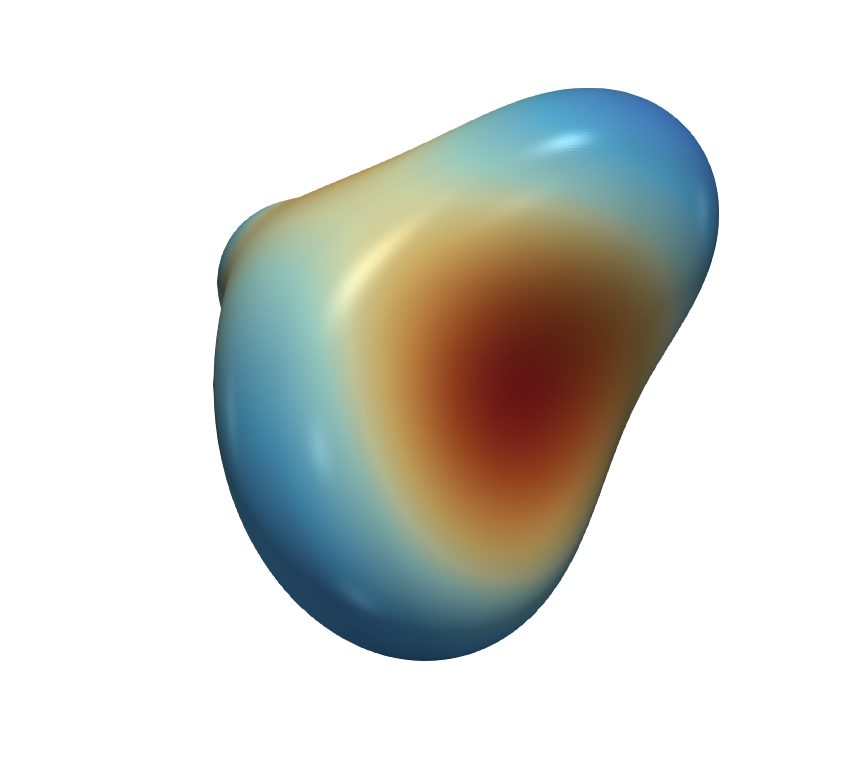}
		\includegraphics{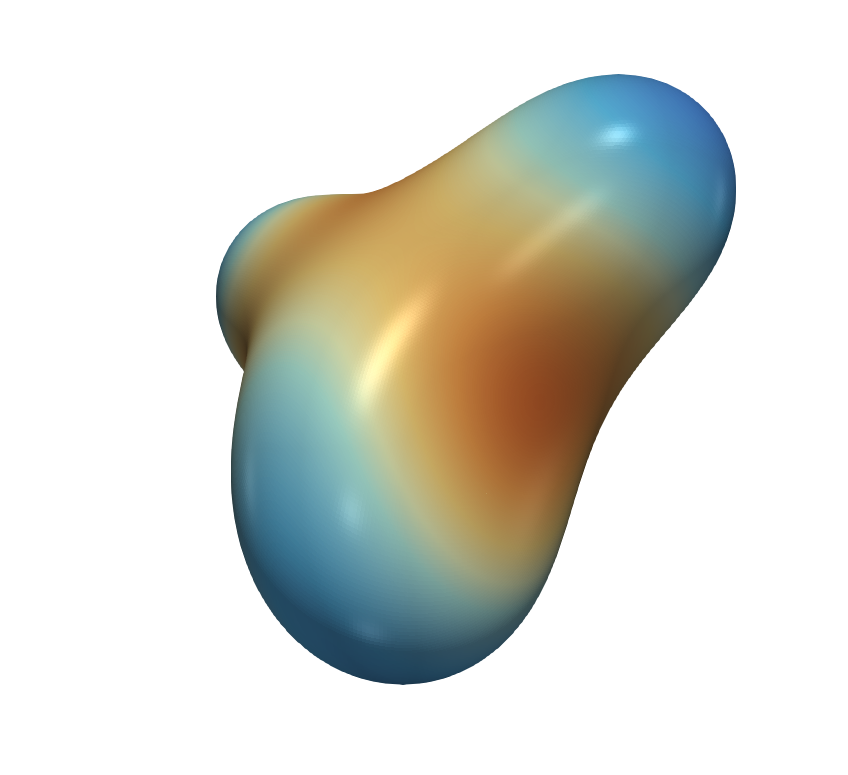}
		\includegraphics{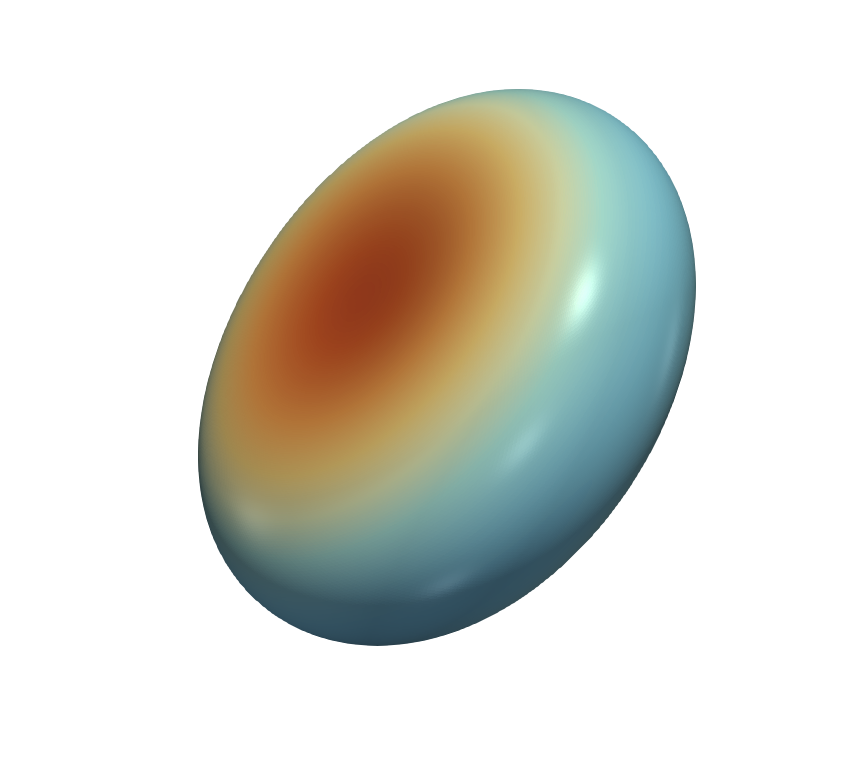}
		\includegraphics{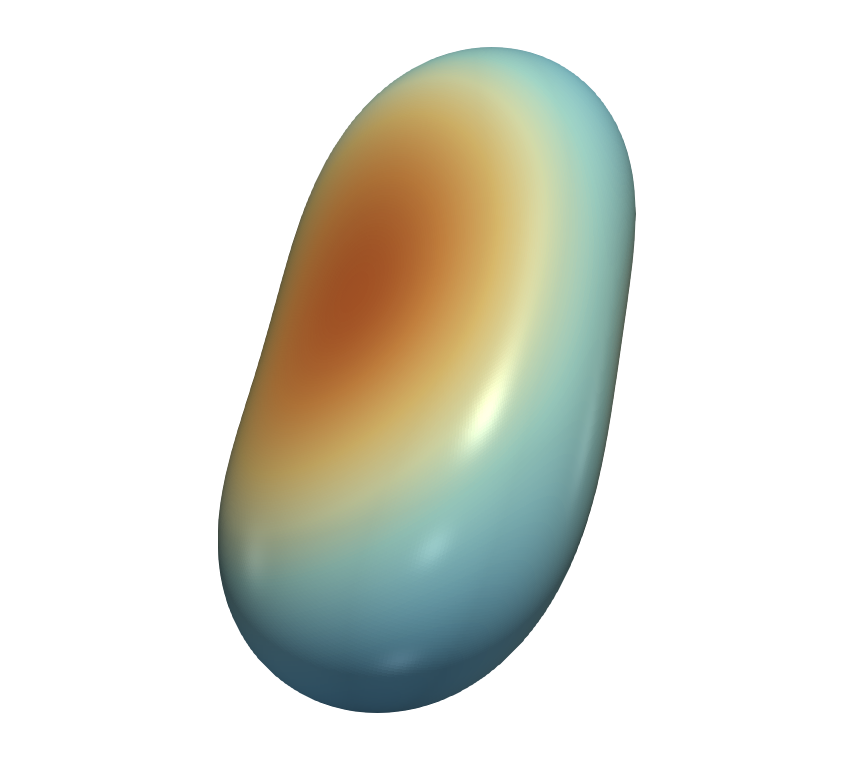}
		\includegraphics{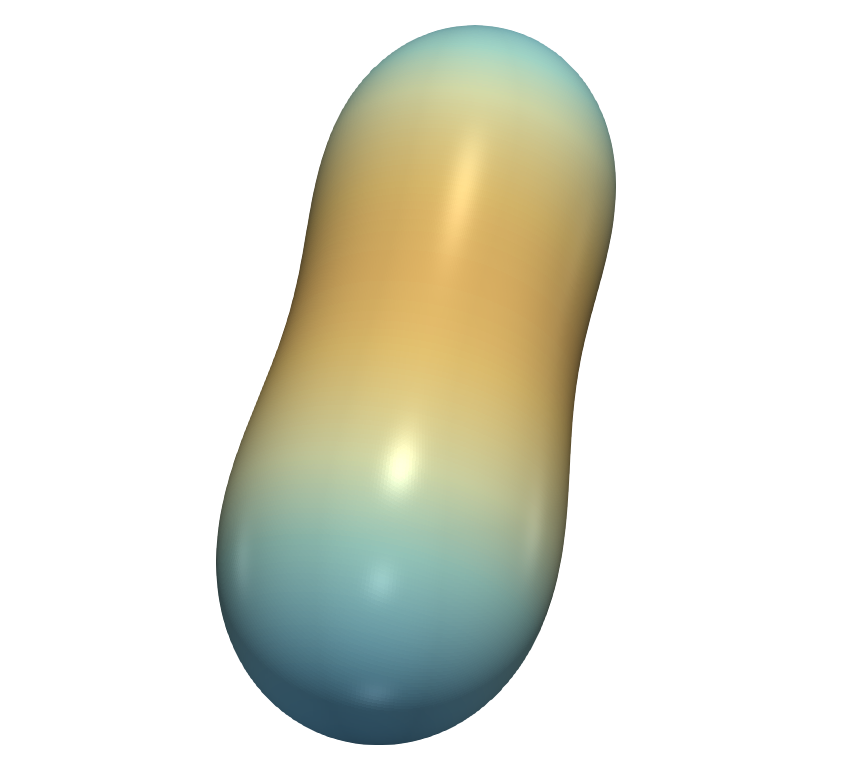}
	\end{adjustbox}

	\begin{adjustbox}{width=0.9\textwidth}
		\includegraphics{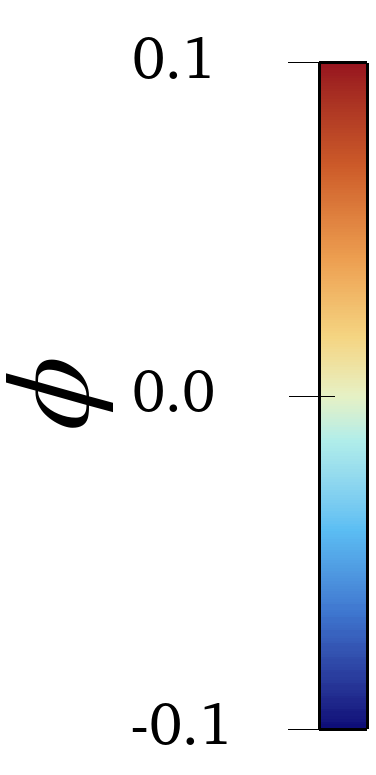}
		\includegraphics{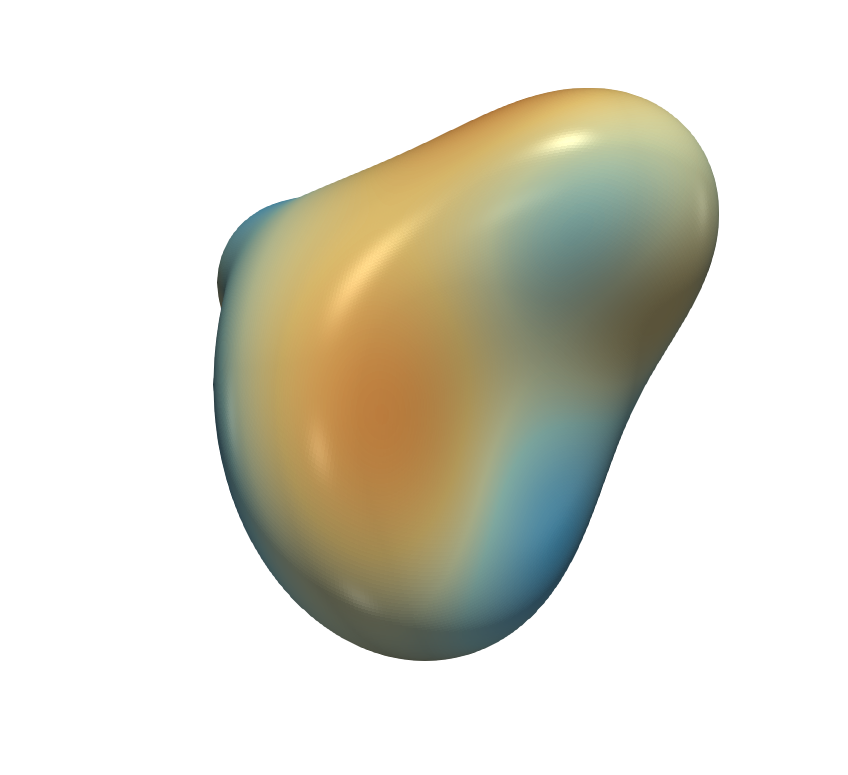}
		\includegraphics{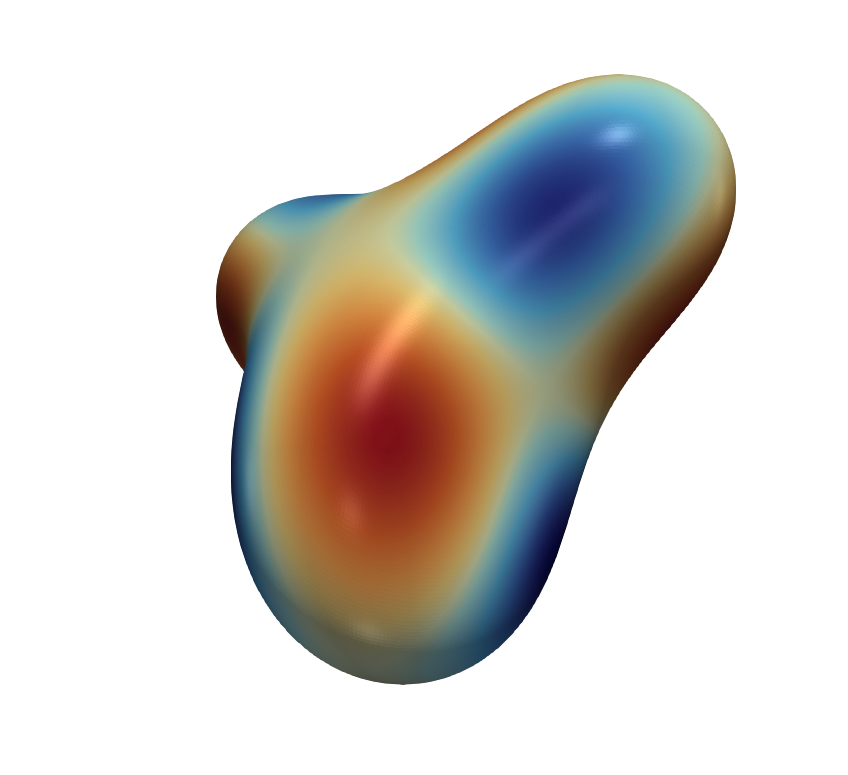}
		\includegraphics{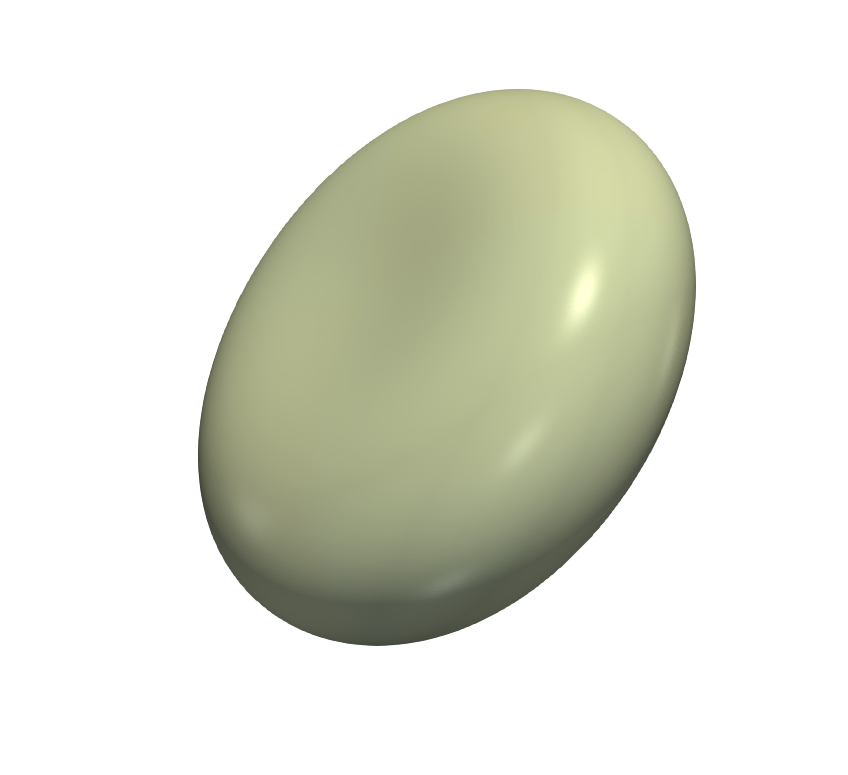}
		\includegraphics{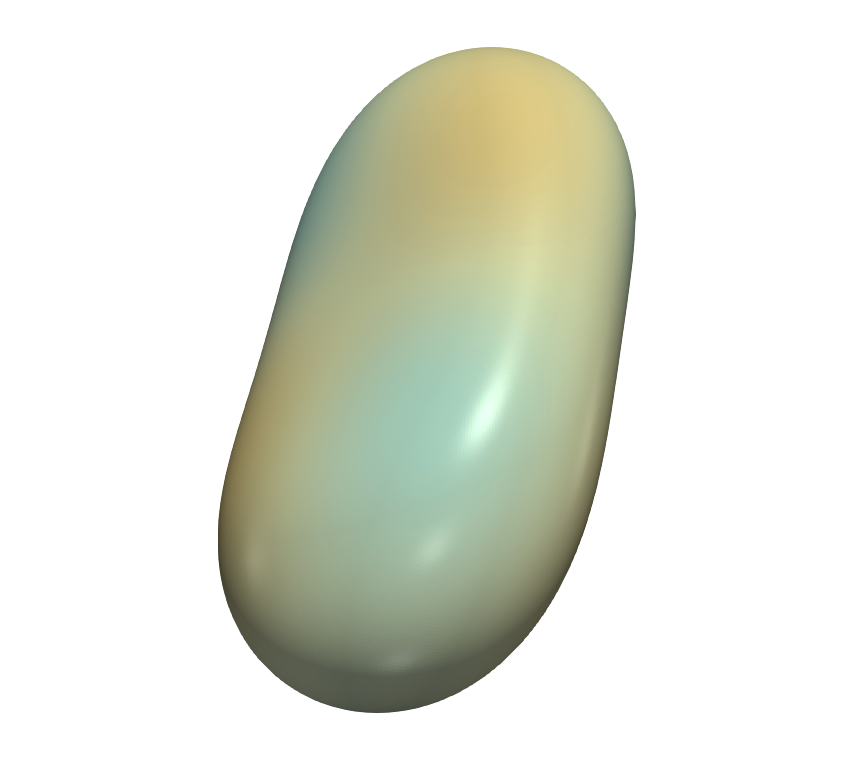}
		\includegraphics{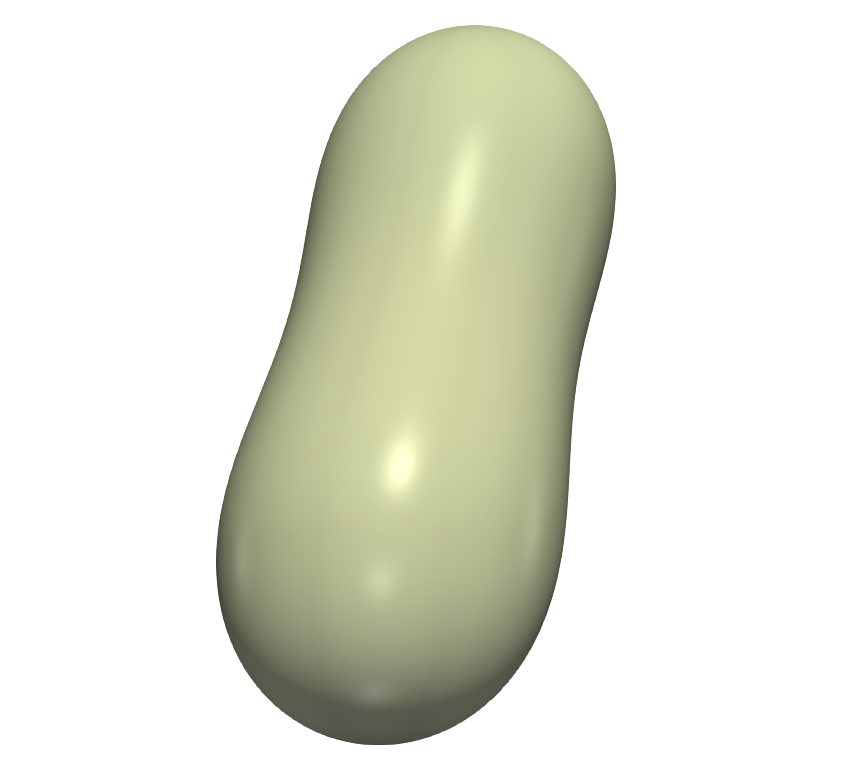}
	\end{adjustbox}

	\begin{adjustbox}{width=0.9\textwidth}
		\includegraphics{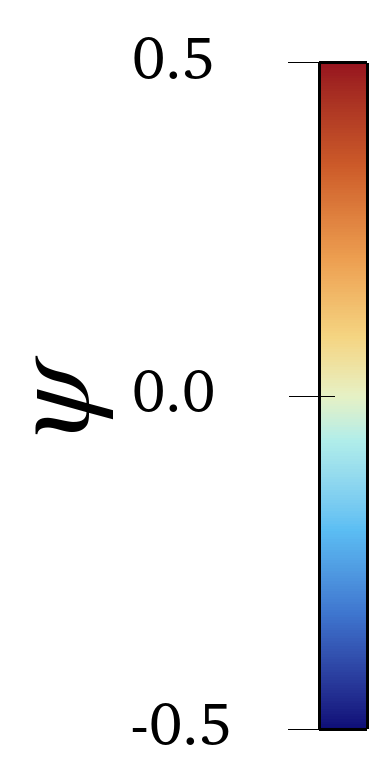}
		\includegraphics{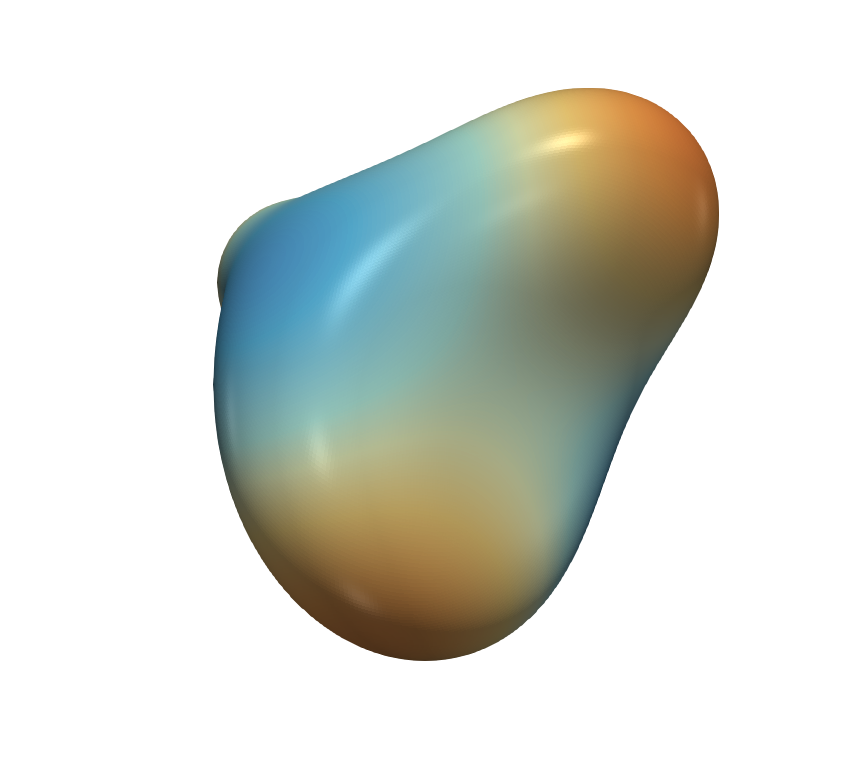}
		\includegraphics{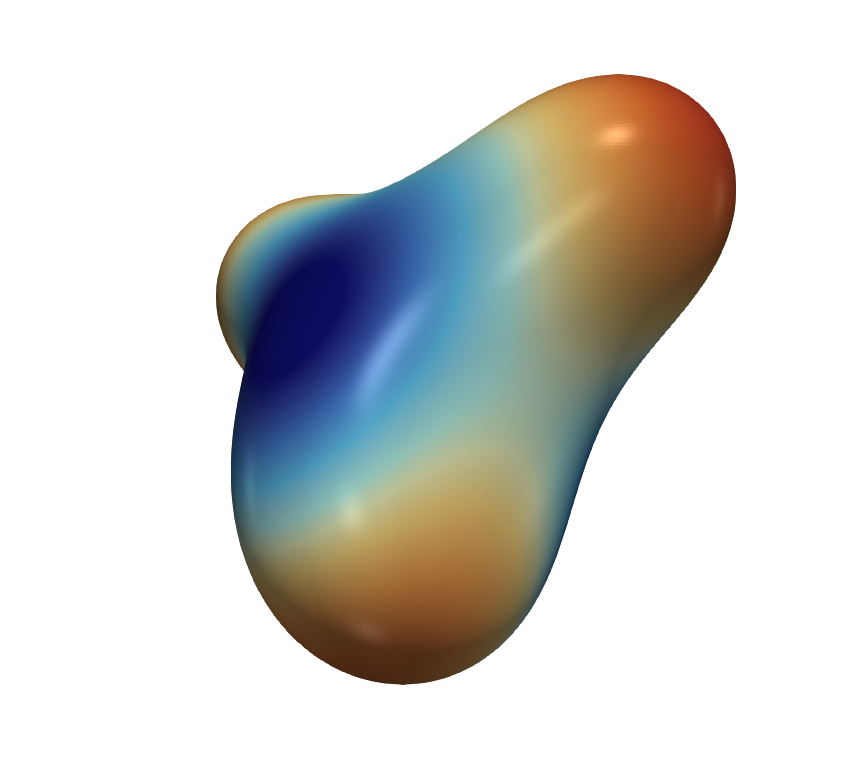}
		\includegraphics{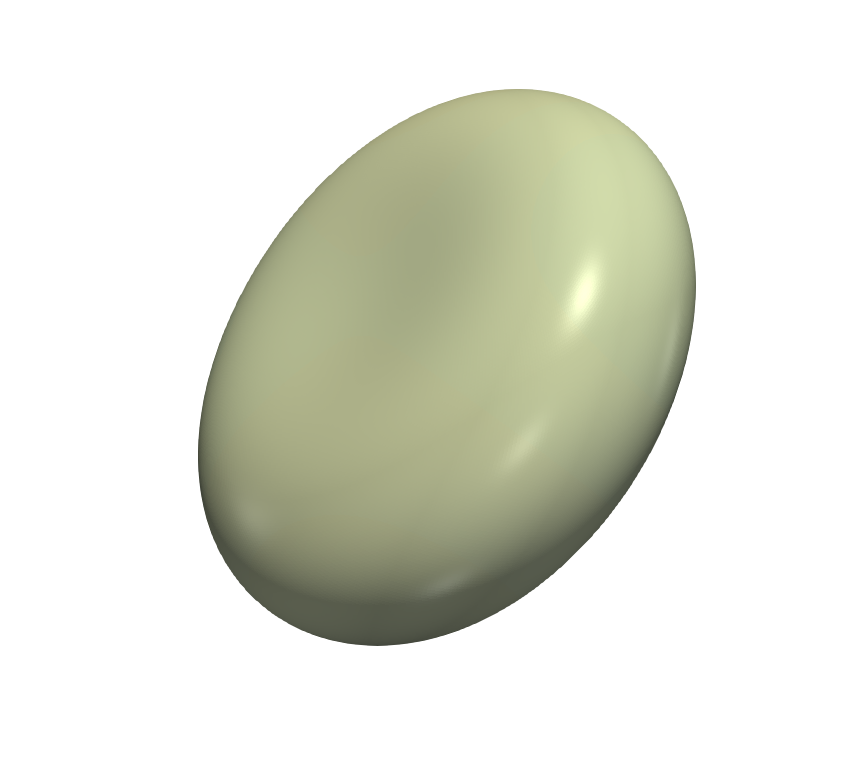}
		\includegraphics{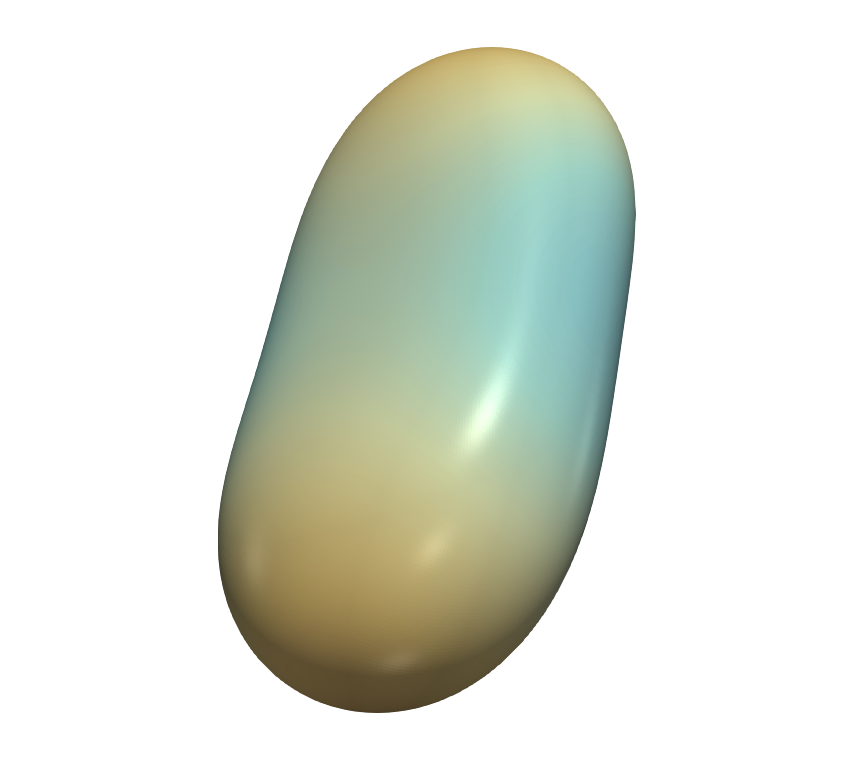}
		\includegraphics{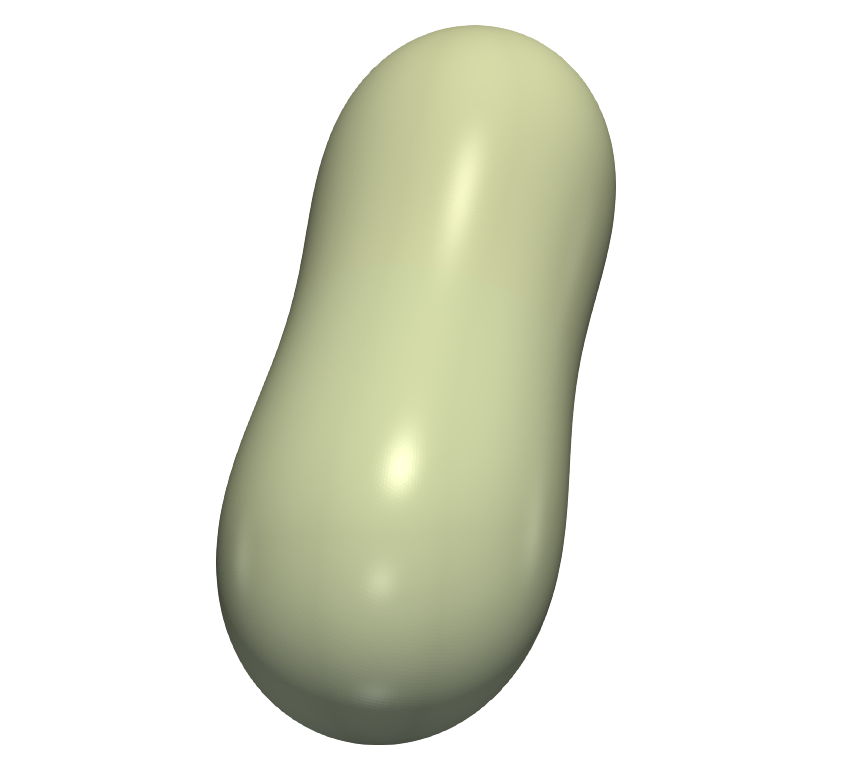}
	\end{adjustbox}

	\begin{adjustbox}{width=0.9\textwidth}
		\includegraphics{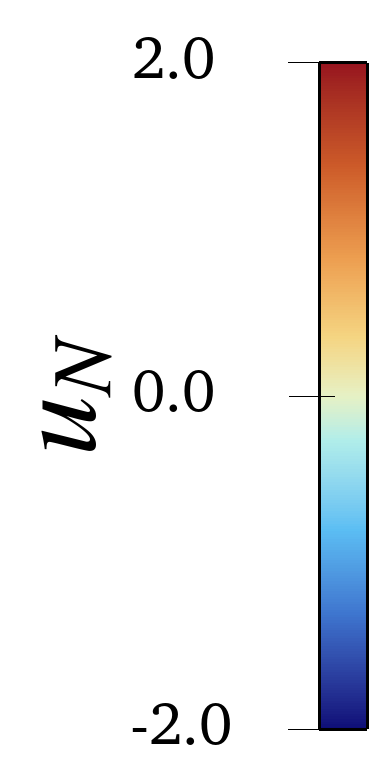}
		\includegraphics{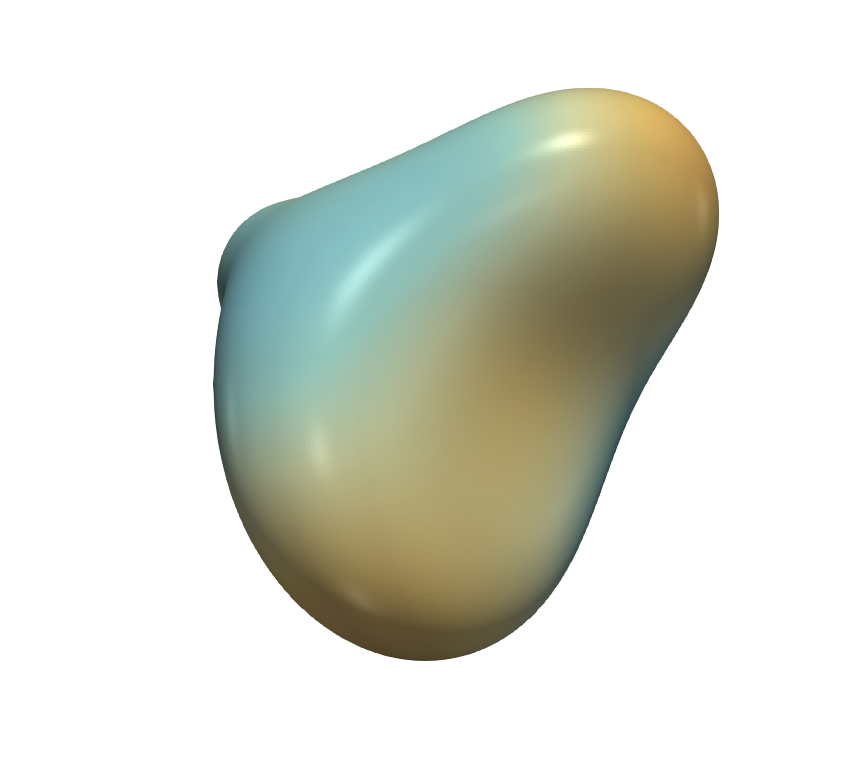}
		\includegraphics{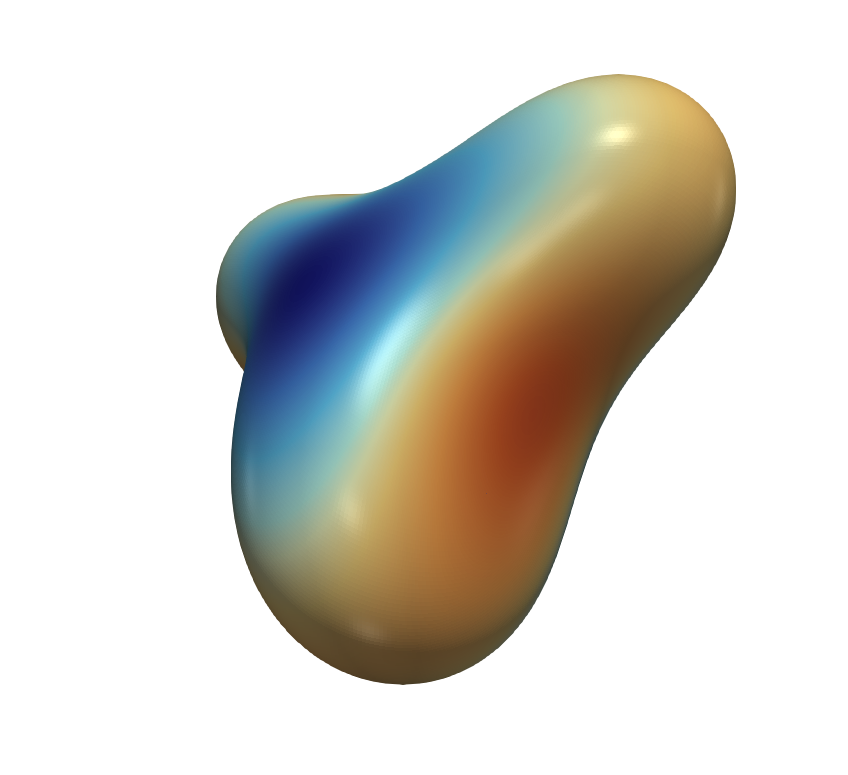}
		\includegraphics{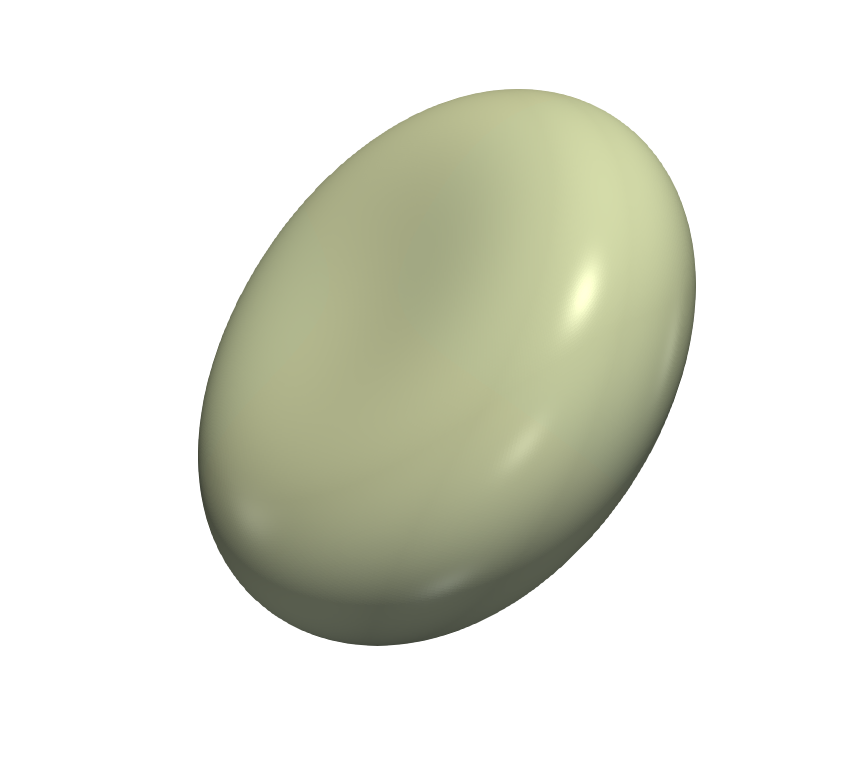}
		\includegraphics{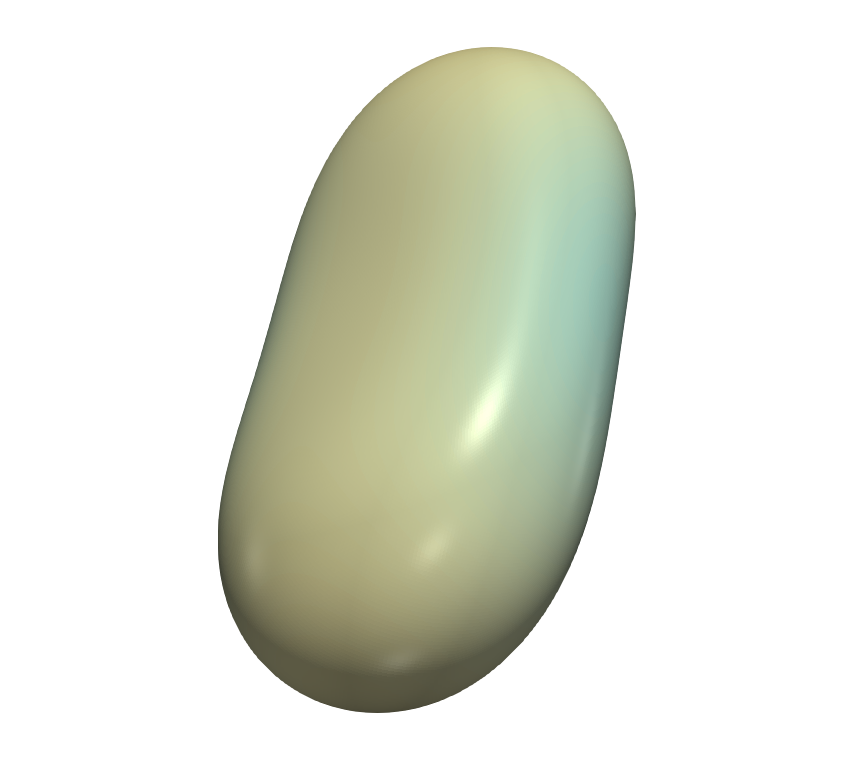}
		\includegraphics{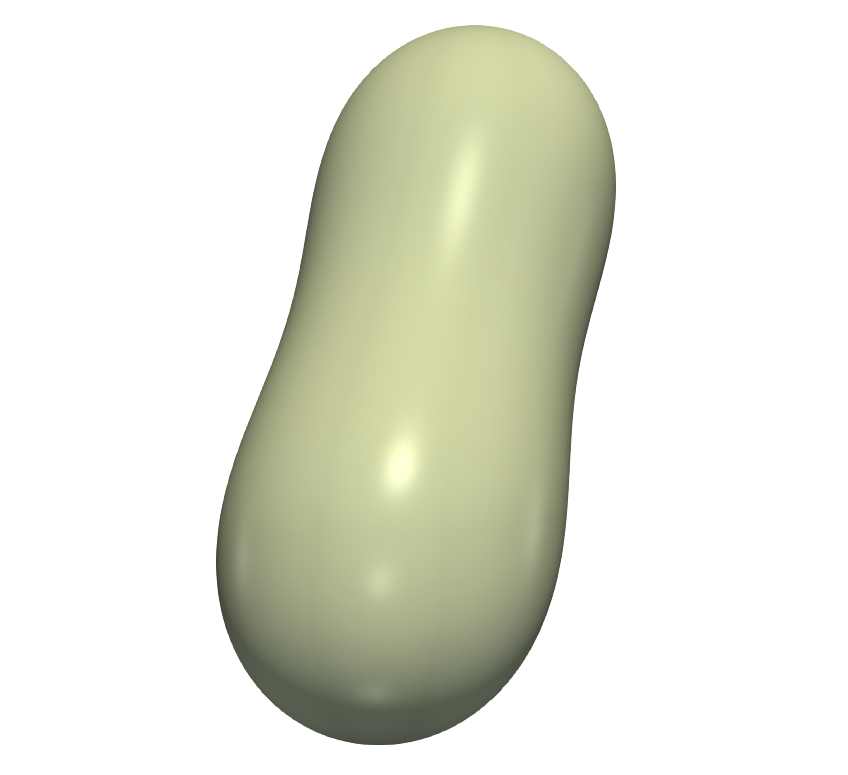}
	\end{adjustbox}

	\begin{adjustbox}{width=0.9\textwidth}
		\includegraphics{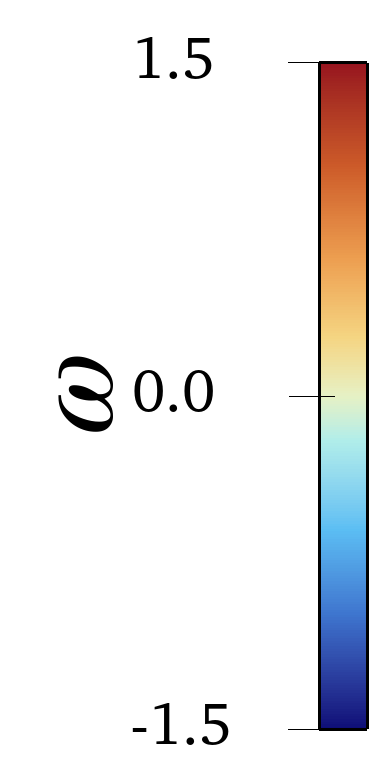}
		\includegraphics{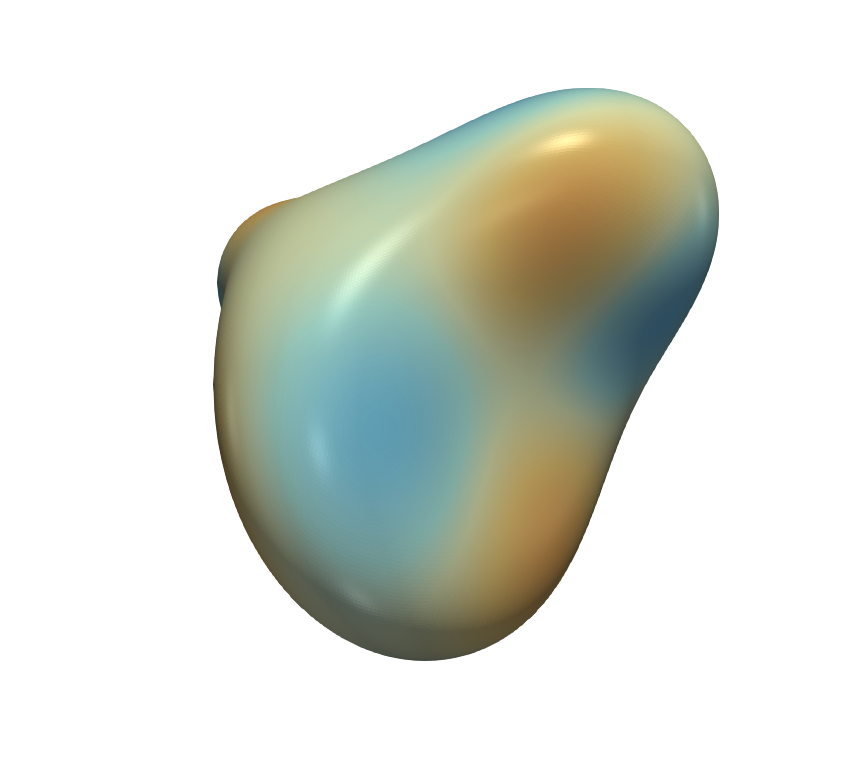}
		\includegraphics{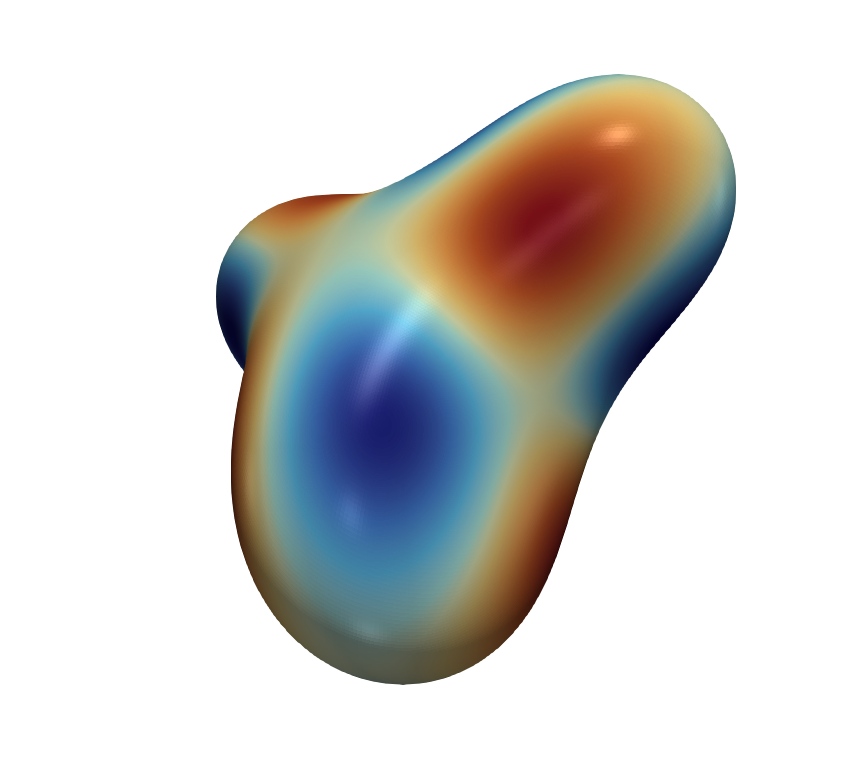}
		\includegraphics{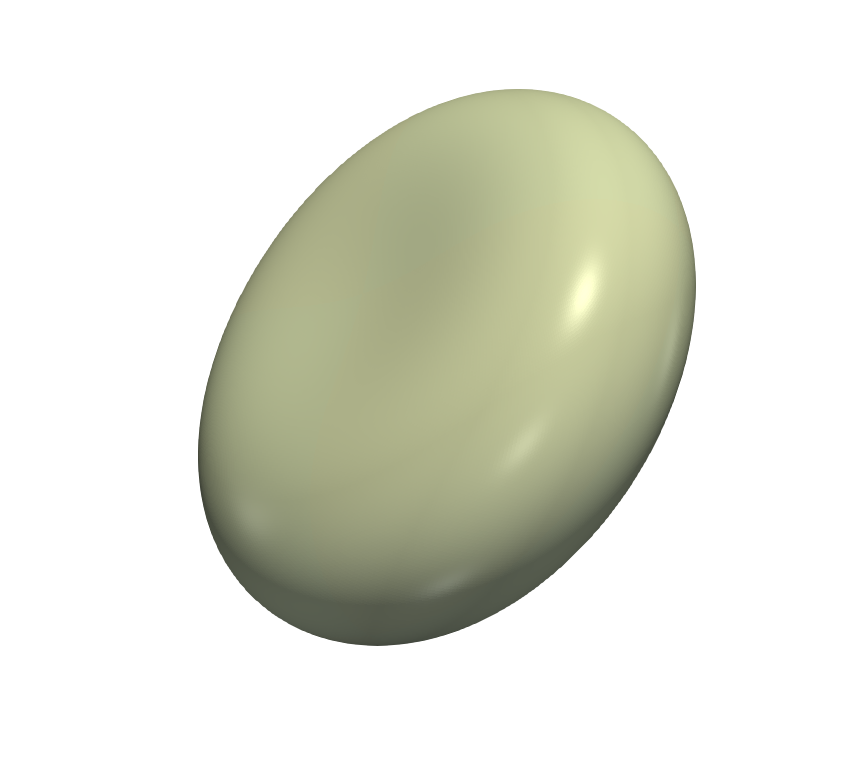}
		\includegraphics{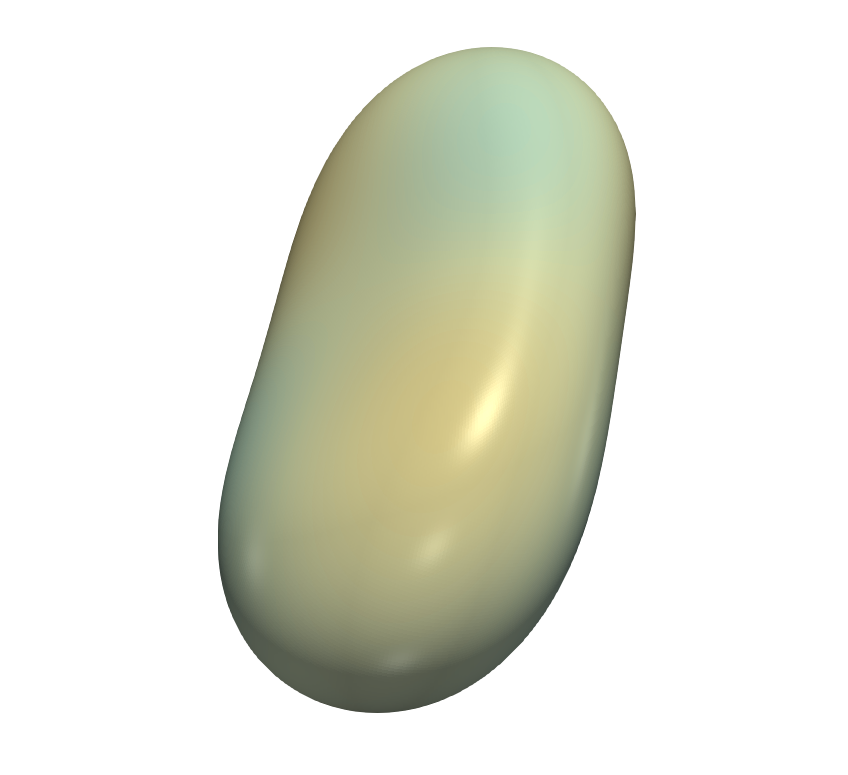}
		\includegraphics{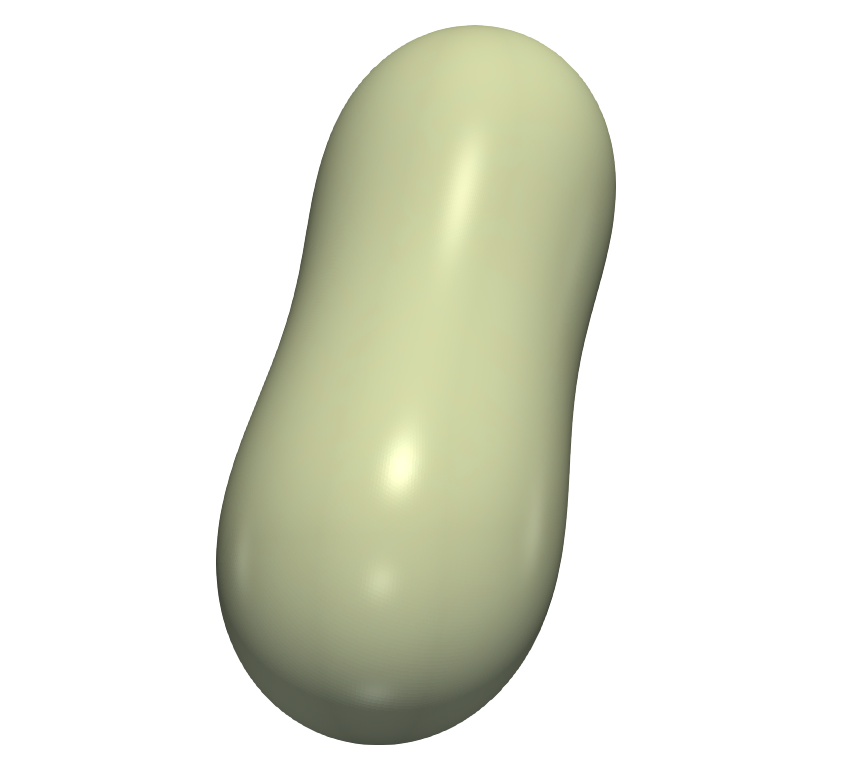}
	\end{adjustbox}

	\begin{adjustbox}{width=0.9\textwidth}
		\includegraphics{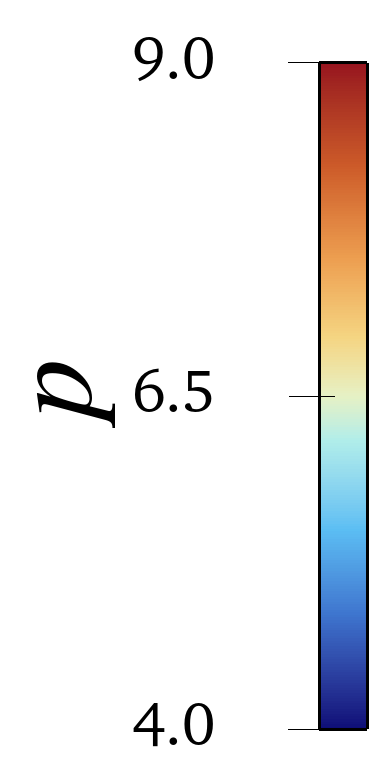}
		\includegraphics{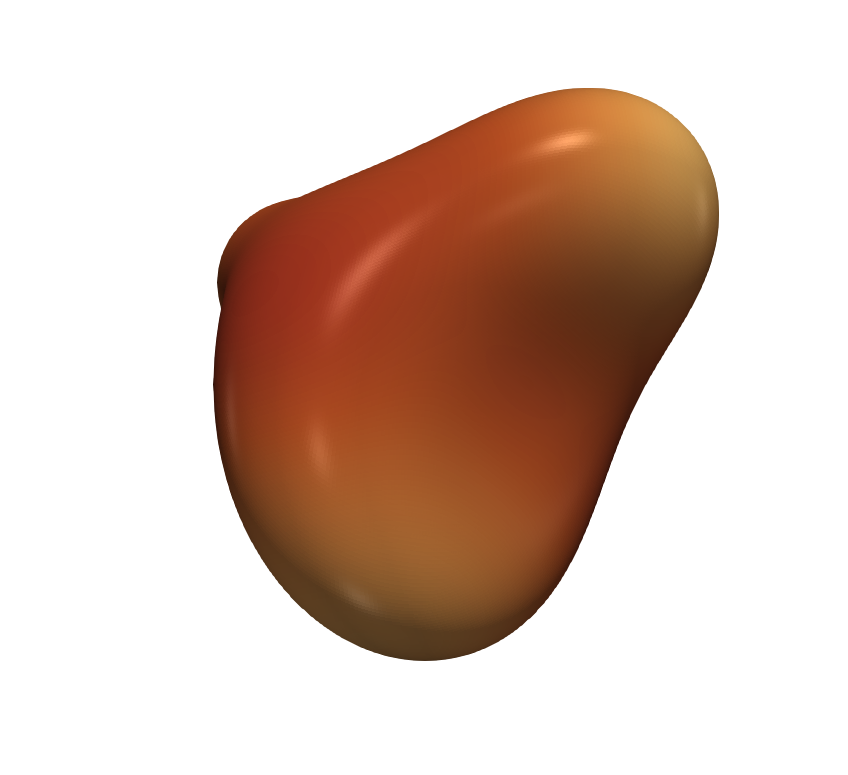}
		\includegraphics{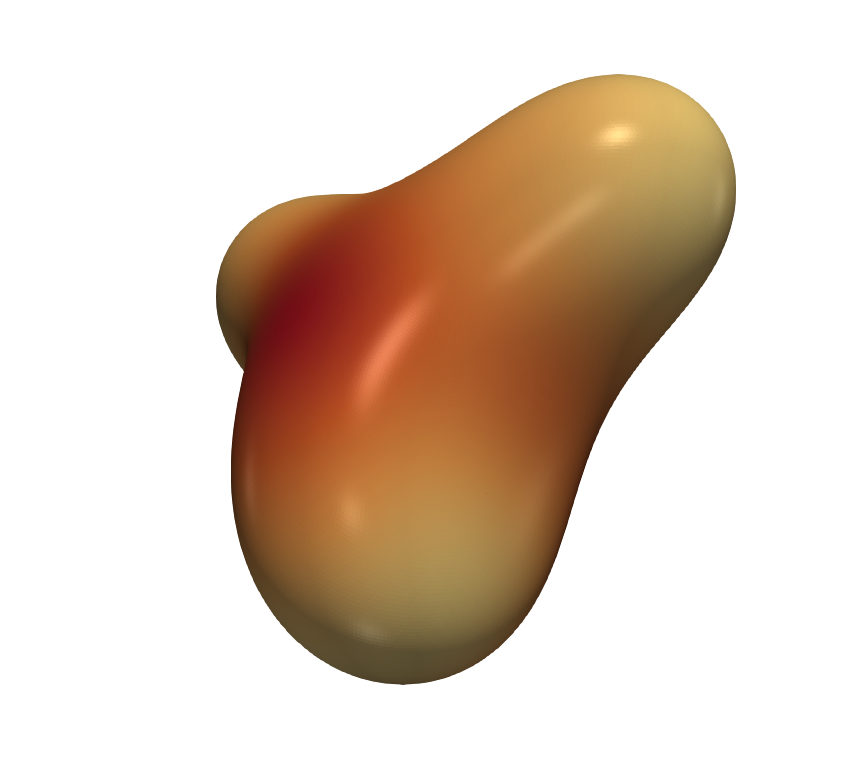}
		\includegraphics{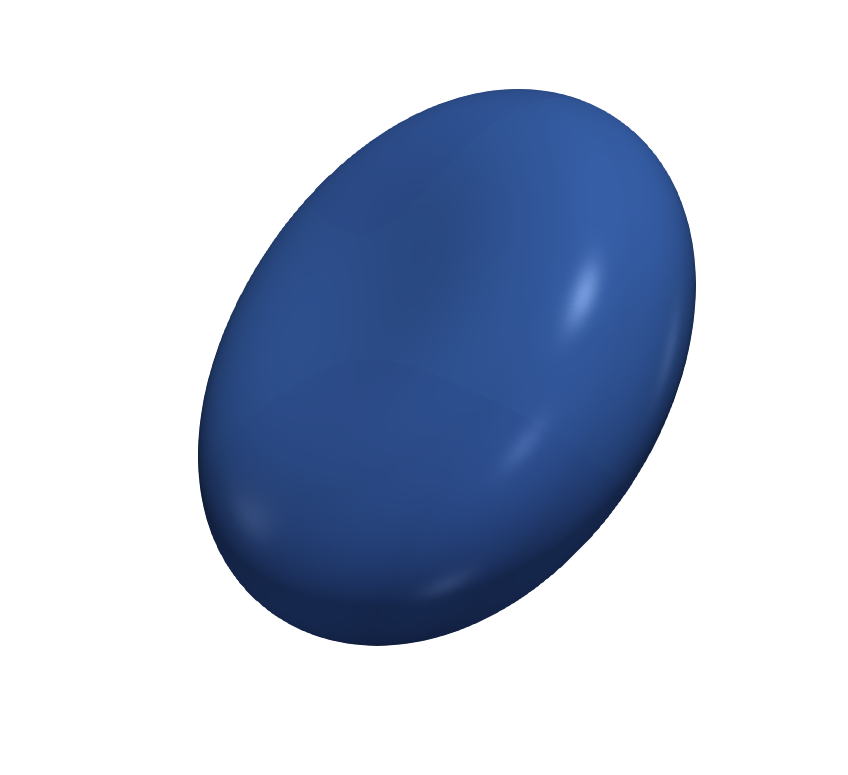}
		\includegraphics{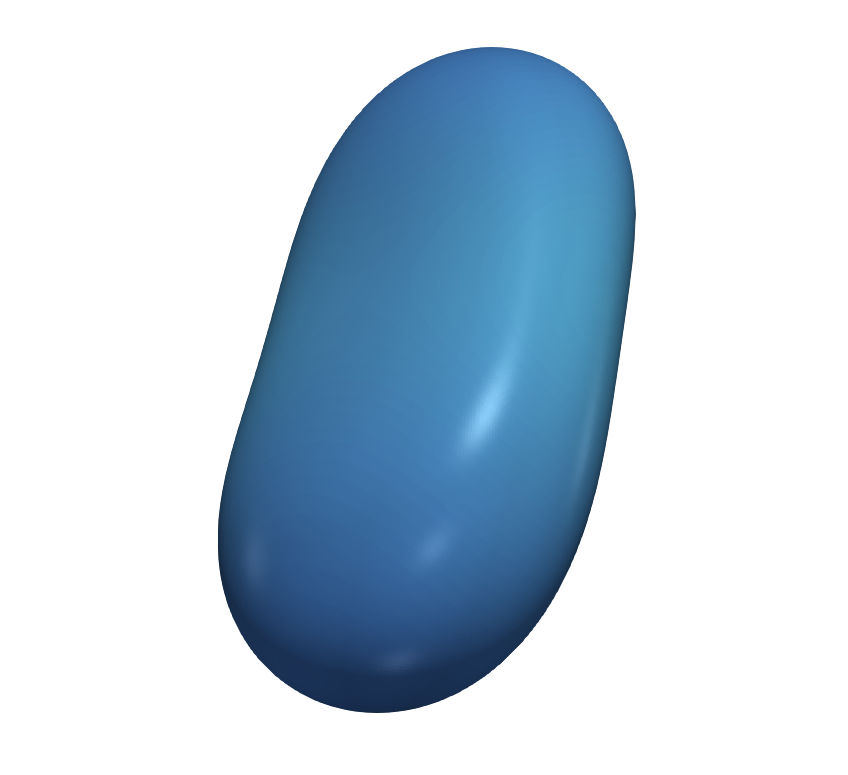}
		\includegraphics{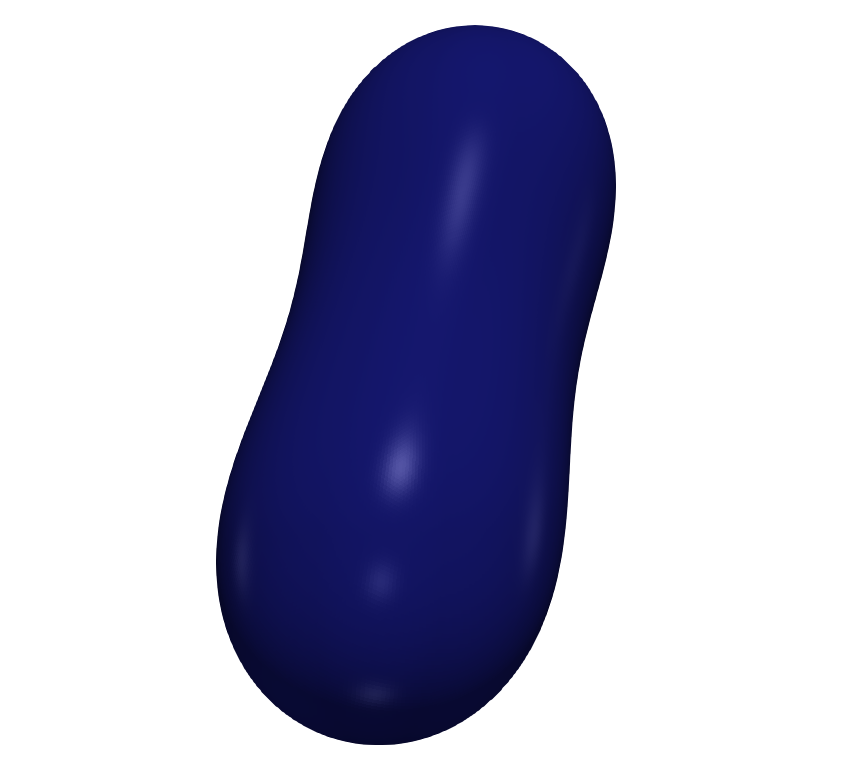}
	\end{adjustbox}

\caption{Time evolution of a perturbed sphere shown as selected time instances with increasing time from left to right, obtained by solving the stream-function formulation of the Surface Stokes-Helfrich model. The time instances $t = 0.09,0.3,3,9;15$ are selected to show characteristic dynamics of relaxation and passing by an oblate shape before reaching the equilibrium prolate shape. Shown is the tangential surface velocity $\flow_{T,h}$ using LIC (line integral convolution) and color coding the magnitude of $\flow_{T,h}$, the mean curvature $\H_h$, the stream-function $\phi_h$, the vector potential $\psi_h$, the normal velocity $u_{N,h}$, the vorticity $\omega_h$ and the pressure $p_h$, from top to bottom. A corresponding movie for $\flow_{T,h}$ is provided in the Electronic Supplement. The numerical parameters are $\tau=0.015$ and $h\approx 0.24$ corresponding to the intermediate values considered in the convergence study in Fig.~\ref{fig:3}.}
\label{fig:1}
\end{figure}

As a benchmark problem, we consider a perturbed sphere parametrization $\param(0)$ with radius depending on spherical coordinates $\theta, \vartheta$ given by $r(\theta, \vartheta) = 1 + r_0\cos \theta \sin 3\vartheta$ with $r_0=\tfrac{2}{5}$. This surface has been used for benchmarking of a Surface Navier-Stokes-Helfrich model previously \cite{krause2023numerical}. We set $\mu=\tfrac{1}{2}$, $\gamma=1$, $\kappa = 1$ and $u_{N,h}^0=0$. The final time $t_{end}$ is chosen to reach the equilibrium configuration. The surface area $A(t) = \int_\S \, \areaelement{}$ and enclosed volume $V(t) = \frac{1}{3} \int_\S \pos \cdot \normal \, \areaelement{}$ correspond to a reduced volume of $\reducedVolume(t) = {6 \sqrt{\pi} \Volume(t)} / {A(t)^{3/2}} \approx 0.81$, which relates to an equilibrium dumbbell shape in the axisymmetric phase diagram, see \cite{seifertShapeTransformationsVesicles1991}. As in
\cite{krause2023numerical} this shape is reached after evolving through the vicinity of an oblate shape, a local minimum of the bending energy $\mathcal{E}_B = \int_\S \frac{\kappa}{2} \H^2 \, \areaelement{}$ with slightly larger bending energy. The non-zero tangential velocity helps to not approach this local minimum and reach the final configuration with zero tangential velocity. Selected time instances of the evolution are shown in Fig.~\ref{fig:1}, showing $\flow_{T,h}$, $\H_h$, $\phi_h$, $\psi_h$, $u_{N,h}$, $\omega_h$ and $p_h$.

\begin{figure}[htpb]
\centering
\begin{adjustbox}{width=\textwidth}
 \includegraphics{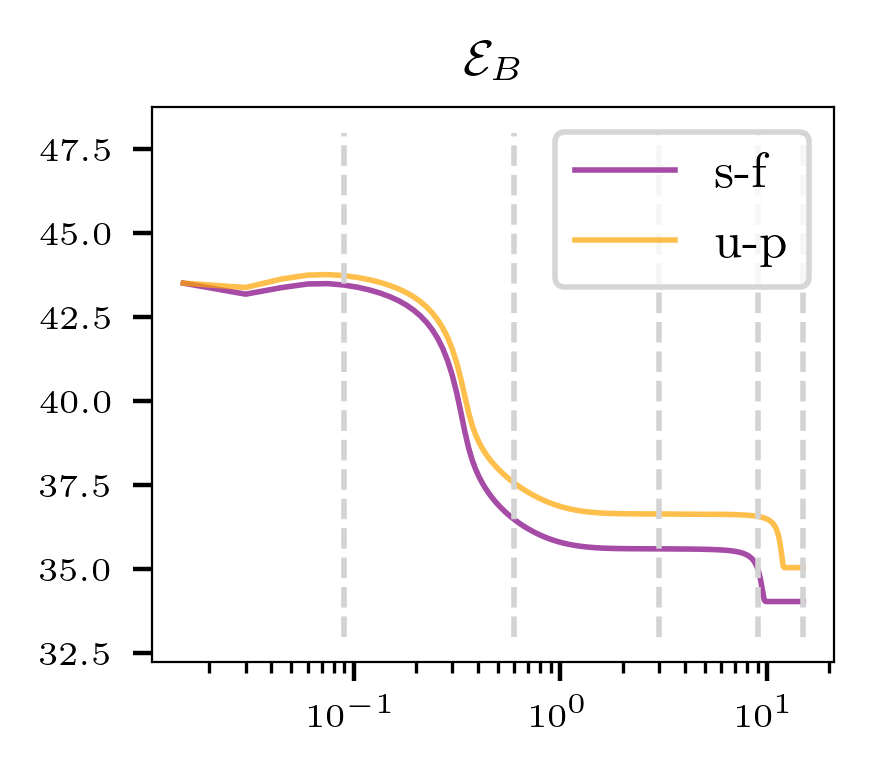}
 \includegraphics{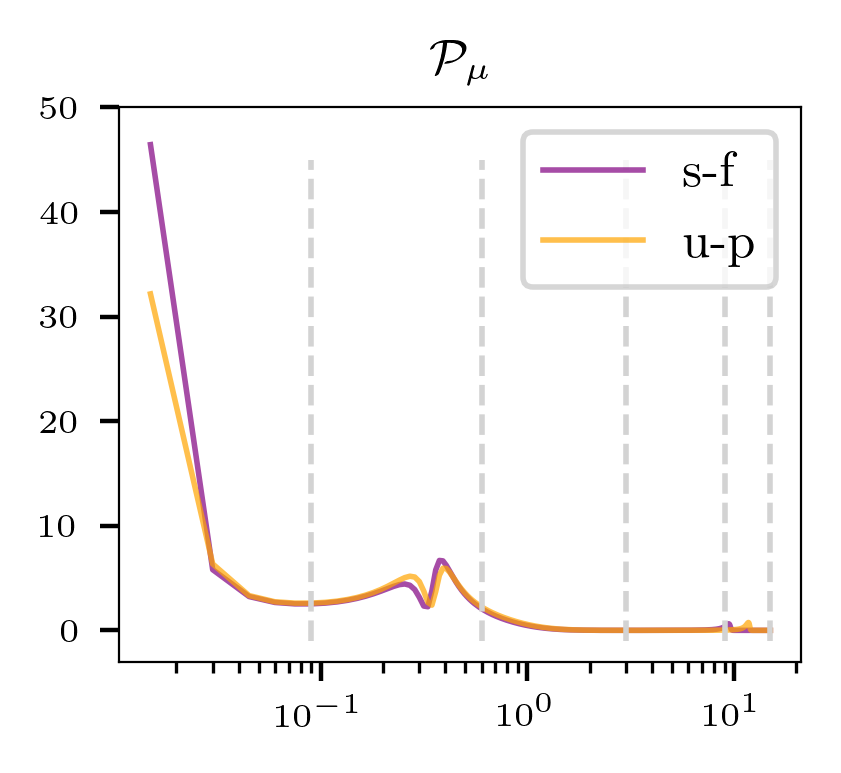}
 \includegraphics{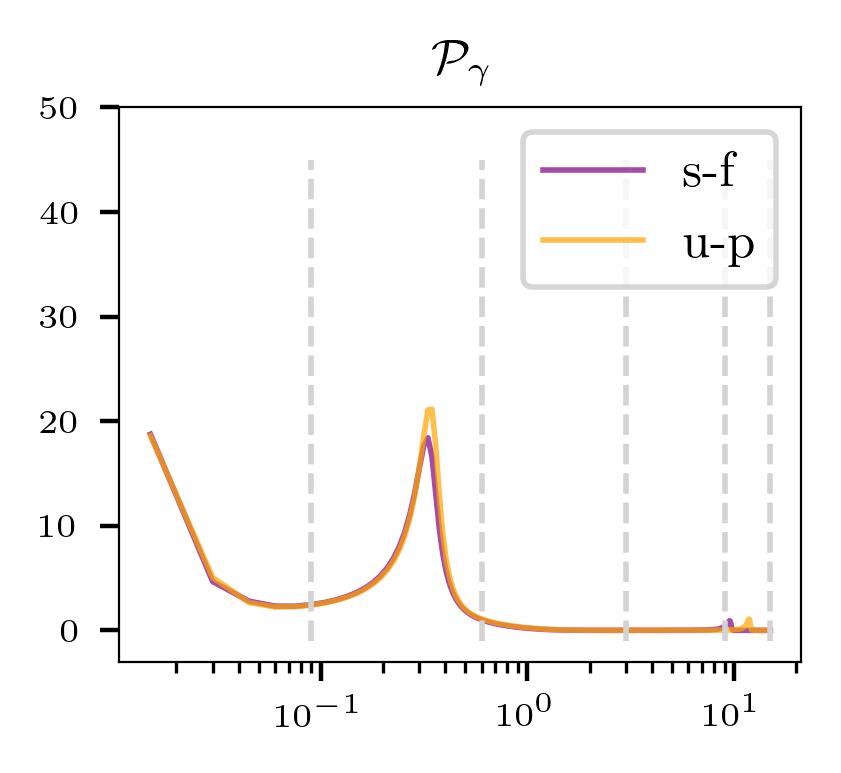}
 \end{adjustbox}
 \begin{adjustbox}{width=\textwidth}
 \includegraphics{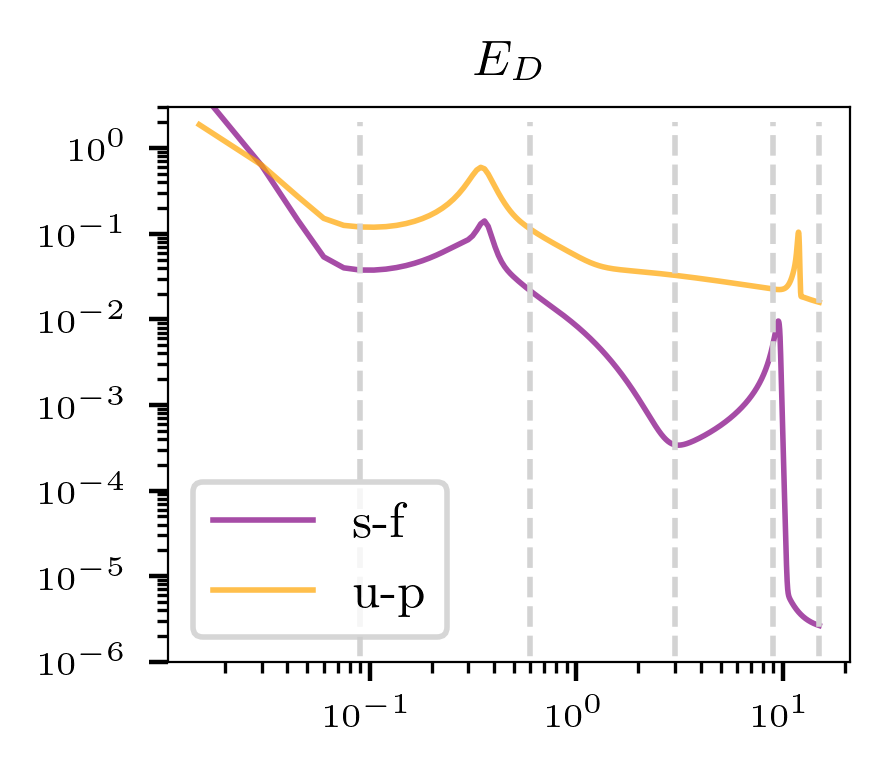}
 \includegraphics{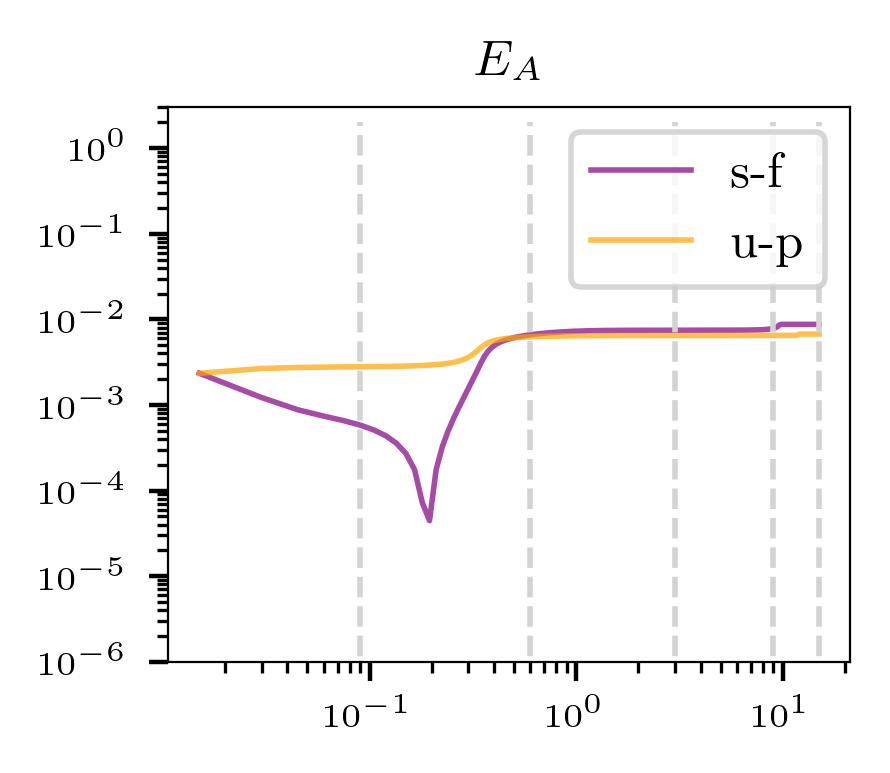}
 \includegraphics{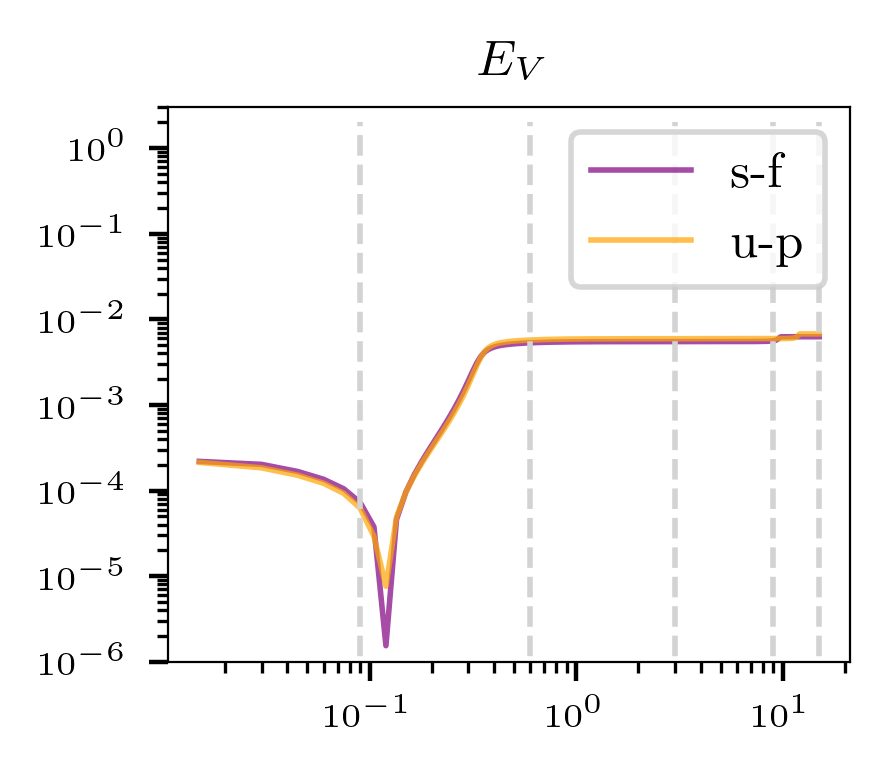}
 \end{adjustbox}
 \caption{Properties of the solution and comparison between the stream-function (s-f) formulation and the velocity-pressure (u-p) formulation of the Surface Stokes-Helfrich model. Shown are the bending energy $\mathcal{E}_B$, the viscous dissipation $\mathcal{P}_\mu$, the frictional dissipation $\mathcal{P}_\gamma$ (top row) and the inextensibility errors $E_D$ and the relative errors in surface area $E_A$ and enclosed volume $E_V$ (bottom row). The time instance considered in Fig.~\ref{fig:1} are marked.}
 \label{fig:2}
\end{figure}

We solve the same problem using the velocity-pressure formulation and compare properties of the solution.  We consider the following properties: the bending energy $\mathcal{E}_B$, the viscous dissipation $\mathcal{P}_\mu = \int_\S \tfrac{\mu}{2} \|\stress(\flow)\|^2 \, d\S $ and the frictional dissipation $\mathcal{P}_\gamma = \int_\S \gamma \|\flow\|^2 \, d\S$. We also consider the inextensibility error $E_D = \| \DivC \flow_h \|_{L^2(\S)}$ as well as relative errors for the conservation of surface area $E_A = |A(t) - A(t_\text{start})|/A(t_\text{start})$ and enclosed volume $E_V = |V(t) - V(t_\text{start})|/V(t_\text{start})$. In Fig.~\ref{fig:2} these quantities are shown over time for both solutions. For the stream-function formulation $\flow_h$ is computed from the reconstructed tangential surface velocity $\flow_{T,h}$, the normal surface velocity $u_{N,h}$ and the normal $\normal_h$, which adds an addition discretization error resulting from the reconstruction.

Both numerical solutions show the desired properties, the bending energy decays and the dissipations more or less coincide. Small deviation are only visible at the transition to the oblate  shape, which is characterized by strong shape deformations and relatively large velocities and the final transition to the prolate shape. The time this transition happens differs between both solutions and also the oblate and prolate shapes slightly differ. These differences are due to different volume/area conservation error behavior, which lead to slightly different reduced volumes $V_r$ at $t_{end}$ and thus different minima of the bending energy. The differences already occur at early stages of the evolution, characterized by strong dissipations. Also the considered inextensibility error shows differences, with larger errors at the beginning of the evolution. Overall the error bounds are comparable between both models and partly even lower for the stream-function approach. To further increase the accuracy, both numerical approaches would benefit from specific structure preserving numerics as recently developed for geometric evolution problems, e.g. Willmore flow \cite{bao2026energystableparametricfiniteelement,Garcke_2026}, which is contained in our problem. However, this goes beyond the topic of this paper, which only aims to establish an alternative numerical approach to the velocity-pressure formulation. For this purpose the reached accuracy is sufficient.

We next demonstrate convergence of the numerical algorithms with respect to refinement of the discretization. As a non-trivial analytical solution of the Surface Stokes-Helfrich model is not known and also numerical analysis results for the algorithms are not yet available, we consider experimental convergence studies towards the numerical solution obtained with the finest resolution. We consider the same problem as before. We define the mesh size $h = \max_T h_T$, where $h_T$ is (ambient) vertex distance of $T$. The initial surface mesh is constructed by global red refinement of an icosahedron and projecting the new vertices to $\param(0)$. The time step $\tau$ is chosen as $h^3$ following previous studies \cite{krause2023numerical,PBV25}. We consider maximum errors over a time interval characterized by strong shape deformations and relatively large velocities and define: $e_D = \| E_D \|_{L^\infty(0.1,1)}$, $e_A = \| E_A \|_{L^\infty(0.1,1)}$ and $e_V = \| E_V \|_{L^\infty(0.1,1)}$. 

\begin{figure}[htpb]
\centering
 \includegraphics[width=0.32\textwidth]{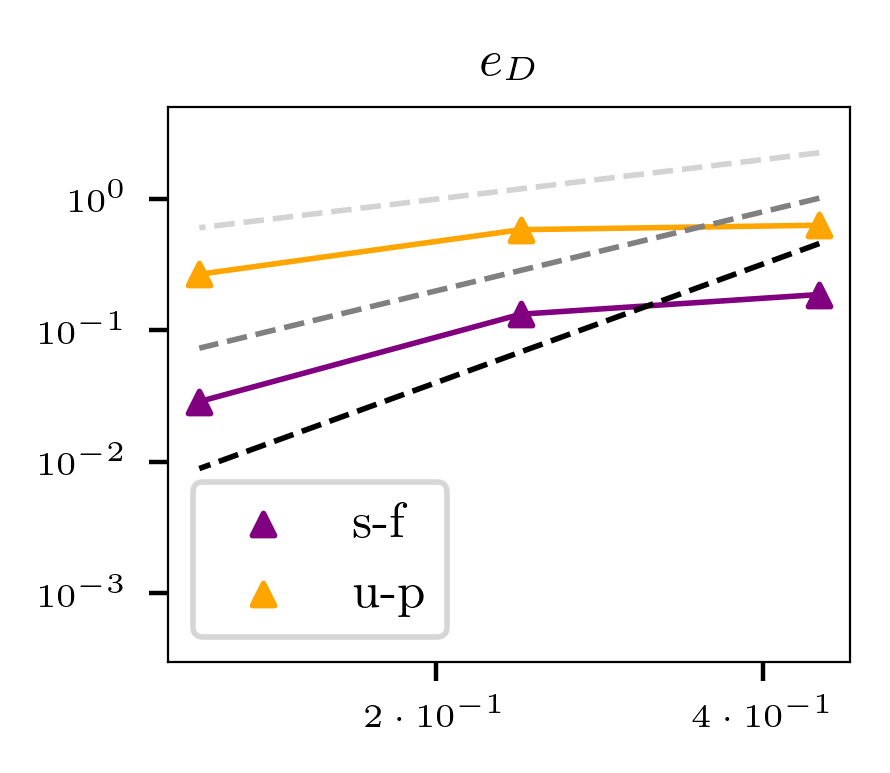}
 \includegraphics[width=0.32\textwidth]{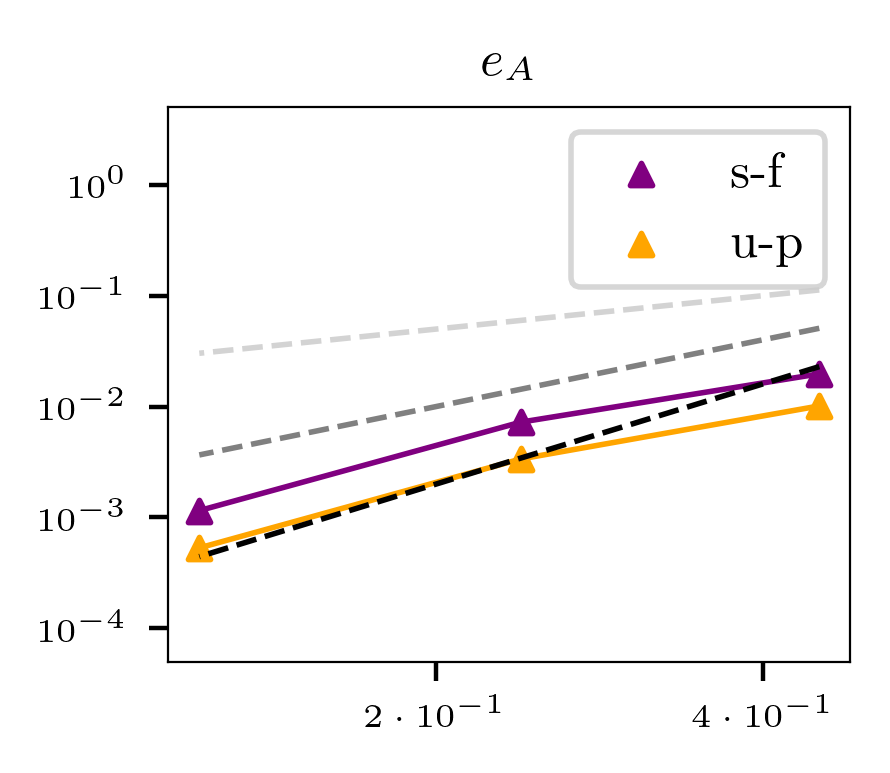}
 \includegraphics[width=0.32\textwidth]{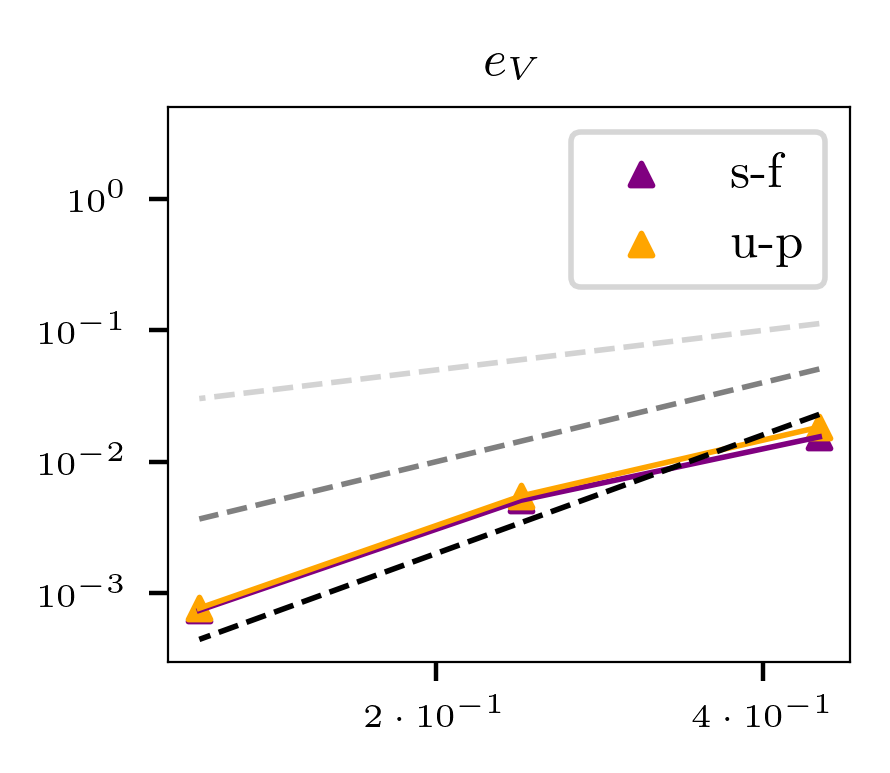}
 \caption{Experimental order of convergence in space and time for $e_D$, $e_A$, $e_V$ 
 for stream-function (s-f) and velocity-pressure (u-p) formulation. Order one, two and three are indicated by dashed lines.}
 \label{fig:3}
\end{figure}

These errors are comparable between both methods with respect to order and absolute values and also agree with previous measures for the velocity-pressure formulation \cite{bachini2023derivation,sischka2025two,PBV25}. Validation against theoretically known orders of convergences are not possible. Only for the surface Stokes problem on a stationary surface results are available for $e_D$, see  \cite{reuskenSfAnalysis,hardering2025parametric,brandner2022finite}. The optimal order for $k=3$ is $3$ for both formulations. This is not achieved for the Surface Stokes-Helfrich problem. The results indicate order $1$ for the velocity-pressure formulation and order $2$ for the stream-function formulation. The reduced order is a consequence of the increased complexity due to various additional unknown $\psi_h, u_{N,h}, p_h, \H_h, \lambda$ and $\mathbf{Y}_h$ and their explicit coupling with geometric quantities through $\K_h, \B_h$ and $\langle \cdot,\cdot\rangle_h$. More detailed explanations would require extending the available numerical analysis results for the surface Stokes problem to the full Surface Stokes-Helfrich problem, which remains an open problem for future research.



\section{Conclusions}
\label{sec:conclusions}

We have proposed a steam-function formulation for the Surface Stokes-Helfrich problem as an alternative to the available velocity-pressure formulation. The formulation is significantly more complex as known stream-function formulations for the Surface Stokes problem on stationary \cite{Nitschke_Voigt_Wensch_2012,10.1093/imanum/dry062,GROSS2018663,torres2019modelling,pearce2019geometrical} or prescribed evolving surfaces \cite{reuther2015interplay,reuther2018erratum}. The main difference is the presence of the surface pressure $p$ in the formulation. We introduced a discretization of the problem using an ALE formulation with standard isoparametric SFEM, a BGN approach to maintain mesh regularity and first order implicit-explicit time-stepping for which we  explored the experimental order of convergence for various quantities and compared the computational results with an established discretization of the velocity-pressure formulation for a non-trivial problem setting.

Based on the obtained errors we conclude that (for the considered test case) the method for the Surface Stokes-Helfrich problem based on the stream-function formulation leads to similar errors and is even slightly more efficient with regard to the inextensibility error than the one based on the velocity-pressure formulation. This is in agreement with the comparison of both methods for the Surface Stokes problem considered in \cite{brandner2022finite}. However, also as in \cite{brandner2022finite}, the difference in efficiency between both methods is not decisive. Comparing the complexity by simply counting the number of unknown (ignoring the lower order for the pressure in the velocity-pressure formulation), the stream-function formulation has one more scalar function to solve for. Additionally considering the reconstruction of the surface velocity further increases the difference. Further comparison would require to consider the computational costs of the sparse direct solver for solving the resulting linear systems. We refrain from such analysis as the results would be highly sensitive to specific settings of the solvers and their ability to deal with the asymmetry in the linear system for the stream-function formulation. We only remark that within the considered settings the computational times for both methods has been very similar.

Overall, the results demonstrate that a stream-function
formulation for the Surface Stokes-Helfrich model and the proposed discretization is feasible. Although additional geometric quantities are required and the pressure cannot be
eliminated, the formulation constitutes a viable alternative to the classical velocity-pressure approach. The main advantages of the stream-function formulation, to circumvent the saddle-point structure and allowing to only solve scalar-valued partial differential equations, remain. While restricted to simply-connected surfaces in the current study, the approach considered in \cite{FluidCohomology,ZhuYinChern2025,bruers2025streamfunction}
for general stationary surfaces can be extended to evolving surfaces by solving also for the harmonic vector field in each time step. However, the technical realization is subject of future research.

\appendix
\section{Derivation of stream-function formulation}
\label{app:derivation}
\newcommand{\generalsca}{f}
\newcommand{\generalscatwo}{g}
\newcommand{\generalvec}{\vec{v}}
Apart from the Helmholtz decomposition \eqref{eq:Helmholtz} for smooth surface tangential vector fields $\flowTest_T$ on simply connected surfaces
\begin{subequations}
\begin{align}
    \flowTest_T = \veccurlS \phi + \gradS \psi \,,\label{eq:Helmholtz_app}
\end{align}
we use the following identities for smooth surface scalar fields $\generalsca$ and $\generalscatwo$
\begin{align}
    \divS(\veccurlS \generalsca) &= 0, \label{id:divcurl}\\
    \curlS(\gradS \generalsca) &= 0, \label{id:curlgrad}\\
    \lapS \generalsca = \divS(\gradS \generalsca) &= \curlS(\veccurlS \generalsca),\label{id:lap} \\
    \veccurlS(\curlS \generalvec_T) &= \vecdivS \left( \P \vecgradC \generalvec_T - (\P \vecgradC \generalvec_T)^T \right) , \label{id:stress1}\\
    \vecdivS (\P \vecgradC \generalvec_T)^T &= \gradS(\divS \generalvec_T) + \K \generalvec_T \label{id:stress2}.
\end{align}
We can write the product rule as
\begin{align}
    \vecdivS(\generalsca \B) &= \B \gradS \generalsca + \generalsca\gradS \H \label{id:prod}
\end{align}
and using the Binet-Cauchy identity we can write
\begin{align}
    \veccurlS \generalsca \cdot \veccurlS \generalscatwo = \gradS \generalsca \cdot \gradS \generalscatwo \label{id:BC}.
\end{align}
From \eqref{id:divcurl} and \eqref{id:curlgrad} it follows
\begin{align}
    \int_\S \veccurlS \generalsca \cdot \gradS \generalscatwo \ \areaelement{} = 0 \label{id:ortho}.
\end{align}
\end{subequations}
The starting point for the proposed stream-function formulation are eqs. \eqref{eq:tmb} - \eqref{eq:incomp}, which we repeat here for convenience. We exclude the volume conservation as this can be reintroduced after reformulation analogously. The equations read:
\begin{subequations}
    \begin{align}
    \gradS p - \mu \left[\vecdivS \rateOfDeformation(\flow_T) - 2\vecdivS(u_N \B)\right] + \gamma \flow_T&= \vec{0}  \label{eq:apptmb} \\
     p\H - \mu \left[\normal \cdot \vecdivC \rateOfDeformation(\flow_T) - 2 \normal \cdot \vecdivC(u_N \B)\right] + \gamma u_N &= f_B \label{eq:appnmb} \\
    \divS \flow_T - u_N \H &= 0 \,.\label{eq:appincomp}
\end{align}
\end{subequations}


Testing \eqref{eq:appincomp} with $\widehat{p}$ and using  \eqref{eq:Helmholtz_app} and \eqref{id:divcurl} to write $\divS \flow_T = \lapS \psi$ we obtain
\begin{equation*}
    \int_\S \lapS\psi\,\widehat p - u_N \mathcal H\, \widehat p \ \areaelement{}
    = \int_\S \gradS\psi\cdot\gradS\widehat p + u_N \mathcal H\,\widehat p \ \areaelement{} = 0\,,
\end{equation*}
and thus
\begin{equation}
    \lapS \psi - u_N \H = 0. \label{eq:sf_deriv_5}
\end{equation}

We can rewrite the covariant divergence of the rate of deformation tensor in \eqref{eq:apptmb} using \eqref{id:stress1}, \eqref{id:stress2} and \eqref{eq:appincomp} as
\begin{subequations}
    \begin{align*}
    \vecdivS \rateOfDeformation(\flow_T) &= \vecdivS \left( \P \vecgradC \flow_T - (\P \vecgradC \flow_T)^T \right) + 2 \vecdivS (\P \vecgradC \flow_T)^T \\
    &= \veccurlS(\curlS \flow_T) + 2 \gradS(\divS \flow_T) + 2\K \flow_T \\
    &= \veccurlS(\curlS \flow_T) + 2 \gradS(u_N \H) + 2\K \flow_T,
\end{align*}
\end{subequations}
testing \eqref{eq:apptmb} with $\widehat{\flow_T}$ and using \eqref{eq:Helmholtz_app} as $\widehat{\flow_T} = \gradS \widehat{\psi} + \veccurlS \widehat{\phi}$ as well as $\flow_T = \gradS \psi + \veccurlS \phi$, we obtain
\begin{align*}
     &\!\int_\S \gradS p \cdot \widehat{\flow_T}\ \areaelement{} -
      \int_\S \mu (\veccurlS(\curlS \flow_T) + 2 \gradS(u_N \H) + 2\K \flow_T) \cdot\widehat{\flow_T}\ \areaelement{} \\
     & \!+\!\int_\S 2\mu \,\vecdivS(u_N \B)\cdot\widehat{\flow_T}\ \areaelement{} + \int_\S \gamma\,\flow_T \cdot \widehat{\flow_T}\ \areaelement{} \\
    = &\!\int_\S \gradS p \cdot \gradS \widehat{\psi} \ \areaelement{}  - \int_\S \mu (2 \gradS(u_N \H) + 2\K \flow_T) \cdot \gradS\widehat{\psi}\ \areaelement{} \\
    & \!-\!\int_\S \mu(\veccurlS(\curlS \flow_T) + 2\K \flow_T) \cdot \veccurlS\widehat{\phi}\ \areaelement{} + \int_\S 2\mu (\B \gradS u_N + u_N \gradS \H ) \cdot \gradS \widehat{\psi} \ \areaelement{} \\
    & \!+\!\int_\S 2\mu (\B \gradS u_N + u_N \gradS \H ) \cdot \veccurlS \widehat{\phi} \ \areaelement{} + \int_\S \gamma \gradS \psi \cdot \gradS \widehat{\psi} \ \areaelement{} + \int_\S \gamma \gradS \phi \cdot \gradS \widehat{\phi} \ \areaelement{} \\
    = &\!\int_\S \gradS p \cdot \gradS \widehat{\psi} \ \areaelement{} - \int_\S 2 \mu \gradS(u_N \H) \cdot \gradS\widehat{\psi} + 2\mu \K \gradS \psi \cdot \gradS\widehat{\psi} + 2\mu\K \veccurlS \phi \cdot \gradS\widehat{\psi}\ \areaelement{} \\
    & \!+\!\int_\S \mu \lapS \phi \cdot \lapS\widehat{\phi} - 2 \mu\K \gradS \phi \cdot \gradS\widehat{\phi} - 2 \mu\K \gradS \psi \cdot \veccurlS\widehat{\phi} \ \areaelement{} \\
     & \!+\!\int_\S 2\mu (\B \gradS u_N + u_N \gradS \H ) \cdot \gradS \widehat{\psi} \ \areaelement{} + \int_\S 2\mu (\B \gradS u_N + u_N \gradS \H ) \cdot \veccurlS \widehat{\phi} \ \areaelement{}  \\
     & \!+\!\int_\S \gamma \gradS \psi \cdot \gradS \widehat{\psi} \ \areaelement{} + \int_\S \gamma \gradS \phi \cdot \gradS \widehat{\phi} \ \areaelement{} \\
     =& \; \bm 0 \,,
\end{align*}
using \eqref{id:ortho}, \eqref{id:prod} and \eqref{id:BC}, as well as \eqref{id:ortho},\eqref{id:BC} and \eqref{id:lap}. We can further test equivalently with $\gradS \widehat{\psi}$ and with $\veccurlS \widehat{\phi}$ leading to
\begin{align*}
    &\int_\S \gradS p \cdot \gradS \widehat{\psi}
    -2 \mu \gradS(u_N \H) \cdot \gradS\widehat{\psi} - 2\mu \K \gradS \psi \cdot \gradS\widehat{\psi} - 2\mu\K \veccurlS \phi \cdot \gradS\widehat{\psi}\ \areaelement{} \nonumber \\
     & +\int_\S 2\mu (\B \gradS u_N + u_N \gradS \H ) \cdot \gradS \widehat{\psi}
    + \gamma \gradS \psi \cdot \gradS \widehat{\psi} \ \areaelement{} \nonumber \\
    = &\; 0 \,,\\
    &\int_\S \mu \lapS \phi \cdot \lapS\widehat{\phi} - 2 \mu\K \gradS \phi \cdot \gradS\widehat{\phi} - 2 \mu\K \gradS \psi \cdot \veccurlS\widehat{\phi} \ \areaelement{} \nonumber\\
    & + \int_\S 2\mu (\B \gradS u_N + u_N \gradS \H ) \cdot \veccurlS \widehat{\phi}
    +  \gamma \gradS \phi \cdot \gradS \widehat{\phi} \ \areaelement{} \nonumber \\
    = & \; 0,
\end{align*}
respectively. Finally, by introducing
\begin{equation}
    \omega = \lapS \phi \label{eq:app:sf4}
\end{equation}
and using \eqref{eq:sf_deriv_5} we obtain the strong formulation
\begin{subequations}
    \begin{align}
    \mu\lapS \omega +2\mu\divS(\K \gradS \phi) + 2\mu \curlS(\K \gradS \psi)
    \nonumber \\
    - 2\mu\curlS(\B \gradS u_N + u_N \gradS \H) - \gamma \omega &= 0 \label{eq:app:sf5a}
    \\ 
    2 \mu\lapS(u_N \H) + 2\mu\divS(\K\gradS \psi) + 2\mu \divS(\K \veccurlS \phi) \nonumber\\
    - 2\mu \divS(\B \gradS u_N + u_N \gradS \H)
    - \lapS p - \gamma u_N \H &= 0 \label{eq:app:sf5b}\\
    \lapS \phi - \omega&= 0. \label{eq:app:sf5c}
\end{align}
\end{subequations}

Testing \eqref{eq:appnmb} with $\widehat{u_N}$ and using \eqref{eq:Helmholtz_app} for $\flow_T = \gradS \psi + \veccurlS \phi$, we obtain
\begin{align*}
     &\int_\S p\H \, \widehat{u_N} - \mu \left[\normal \cdot \vecdivC \rateOfDeformation(\flow_T) - 2 \normal \cdot \vecdivC(u_N \B)\right] \, \widehat{u_N} + \gamma u_N \, \widehat{u_N} \ \areaelement{} \\
    = & \int_\S p\H \, \widehat{u_N} + \mu \rateOfDeformation(\flow_T) : \gradC ( \widehat{u_N} \normal) - 2 \mu(u_N \B) : \gradC (\widehat{u_N} \normal) + \gamma u_N \, \widehat{u_N} \ \areaelement{} \\
    = & \int_\S p\H \, \widehat{u_N} - \mu \rateOfDeformation(\flow_T) : \B \widehat{u_N} + 2 \mu u_N \B : \B \widehat{u_N} + \gamma u_N \, \widehat{u_N} \ \areaelement{} \\
    = & \int_\S p\H \, \widehat{u_N} - 2\mu (\P \vecgradC \flow_T) : \B \widehat{u_N} + 2 \mu u_N \B : \B \widehat{u_N} + \gamma u_N \, \widehat{u_N} \ \areaelement{} \\
    = & \int_\S p\H \, \widehat{u_N} + 2\mu \flow_T \cdot \vecdivS(\B \widehat{u_N}) + 2 \mu u_N \B : \B \widehat{u_N} + \gamma u_N \, \widehat{u_N} \ \areaelement{} \\
    = & \int_\S p\H \, \widehat{u_N} + 2\mu (\gradS \psi + \veccurlS \phi) \cdot (\B \gradS \widehat{u_N} + \widehat{u_N}\gradS \H) \ \areaelement{} \\
    & + \int_\S 2 \mu u_N \B : \B \widehat{u_N} + \gamma u_N \, \widehat{u_N} \ \areaelement{} \\
    = & \int_\S f_B \, \widehat{u_N} \ \areaelement{}\,,
\end{align*}
where we have used the symmetry of the shape operator as well as \eqref{id:prod}. We thus obtain the strong formulation
\begin{align}
    -2 \mu \divS(\B (\gradS \psi + \veccurlS \phi)) + 2 \mu (\gradS \psi + \veccurlS \phi) \cdot \gradS \H &\nonumber \\
    +2 \mu u_N \B : \B + p\H + \gamma u_N +\lambda - f_B  &= 0 \,. \label{eq:app:sf6}
\end{align}

Summarizing \eqref{eq:app:sf5a}, \eqref{eq:app:sf5b}, \eqref{eq:app:sf6}, \eqref{eq:app:sf5c}, \eqref{eq:app:sf4} and \eqref{eq:sf_deriv_5} represent the equations of the stream-function formulation in Problem \ref{problem:sf} (neglecting the volume constraint).

\section*{Acknowledgments}
 A.V. acknowledges support from the German Research Foundation (DFG), project numbers 417223351 (FOR3013). The authors gratefully acknowledge the computing time made available to them on the high-performance computer at the NHR Center of TU Dresden. This center is jointly supported by the Federal Ministry of Research, Technology and Space of Germany and the state governments participating in the \href{www.nhr-verein.de/unsere-partner}{NHR}.

\bibliographystyle{siamplain}
\bibliography{references}
\end{document}